\documentclass{article}

\usepackage{amsfonts}
\usepackage{amsmath}
\usepackage{enumerate}
\usepackage{epsfig}
\usepackage{amssymb}
\usepackage[T1]{fontenc}
\usepackage{bbm}

\newtheorem{theorem}{Theorem}
\newtheorem{lemma}[theorem]{Lemma}%

\newtheorem{remark}{Remark}%

\newtheorem{definition}{Definition}%

\newcommand{\one}{\mathbbm{1}}
\begin{document}

\title{Multiple-Population Discrete-Time Mean Field Games with Discounted and Total Payoffs: The Existence of Equilibria}


\author{Piotr Wi\k{e}cek\\Faculty of Pure and Applied Mathematics,\\Wroc{\l}aw University of Science and Technology,\\ Wybrze{\.z}e Wyspia{\'n}skiego 27, 50-370 Wroc{\l}aw, Poland\\
email: Piotr.Wiecek@pwr.edu.pl}
\maketitle
%
%

%
%


\abstract{In the paper we present a model of discrete-time mean-field game with several populations of players. 
Mean-field games with multiple populations of the players have only been studied in the literature in the continuous-time setting. The main results of this article are the first stationary and Markov mean-field equilibrium existence theorems for discrete-time mean-field games of this type. 
We consider two payoff criteria: $\beta$-discounted payoff and total payoff. The results are provided 
under some rather general assumptions on one-step reward functions and individual transition kernels of the players. In addition, the results for total payoff case, when applied to a single population, extend the theory of mean-field games also by relaxing some strong assumptions used in the existing literature.}

\maketitle

\section{Introduction.}
\label{intro}

Dynamic games with a large number of players are a natural tool to model dynamic interactions in many areas of science, yet they do not attract much attention due to complexity of such models. One of the natural ways to deal with problems with a large number of agents that have been developped in different fields of research is to replace such complex models with relatively simpler ones with a continuum of infinitesimal players. This kind of approximations have appeared in one-step games at least since two seminal papers by Wardrop \cite{Wardrop} and Schmeidler \cite{Schmeidler}, but for a long time have not been introduced to dynamic game models.
The situation has changed since a series of papers by Lasry and Lions \cite{LL2007a,LL2007b} and by Huang, Caines and Malham\'{e}
\cite{HCM2006,HCM2007a,HCM2007c} where models of non-cooperative differential games with a continuum of identical players have been introduced. The idea on which these models were founded was that for games of this type in the limit with infinite number of players, the game problem can be reduced to a much simpler single-agent decision problem. A huge number of publications on the topic have followed during the last decade and the literature is still growing fast. A review of the existing results on differential-type mean-field games can be found in the books \cite{BFY,CarmonaDelarue} or the survey \cite{Gomes_survey}.

Similar discrete-time models have appeared in the literature significantly earlier in the paper by Jovanovic and Rosenthal \cite{JovanovicRosenthal} under the name of anonymous sequential games, but have not attracted as much attention as their continuous-time counterparts. However, since then some further theoretical results on this type of games have appeared. The models with discounted payoff criterion have been studied in \cite{BB92,BB95,Chakrabarti03,AdlakhaJohari,SaldiBasarRaginsky,ElliotLiNi,SaldiBasarRaginsky18,SaldiBasarRaginsky19,SaldiBasarRaginsky22}. Conditions under which Nash equilibria in finite-player discounted-utility games converge to equilibria of respective anonymous models were analyzed in \cite{Green80,Green84,Housman88,Sabourian90,SaldiBasarRaginsky}. In \cite{PiWEA,PiW_new,Saldi} long-time average payoff has been considered, while \cite{PiWEA} have treated the games with total reward criterion. In \cite{AnahtarciKariksizSaldi,AnahtarciKariksizSaldi22}, algorithms allowing to compute mean-field equilibria in both discounted and average reward games have been presented.

All of the papers enumerated above have considered the case with only one population of symmetrical players. There is no reason however not to consider mean-field games with a bigger number of populations. As long as this number is small, considering this kind of limit model rather than a game with a huge finite number of players should be a significant simplification of the problem. In case of continuous time, this type of models have been introduced in \cite{HCM2006} and further studied in \cite{Feleqi,Cirant2015,Cirant2017,AchdouBardi,BardiCirant,BensoussanHuangLauriere,MorandottiSolombrino}. As far as we know, there have been no papers on discrete-time mean-field games with multiple populations of players. In this article, we try to fill in that gap by introducing two models of games of this type: one with discounted payoff, another with total payoff. In both cases we provide the results 
about the existence of mean-field equilibria in such games under some natural assumptions. 
It is worth mentioning here that some of the results we present, notably all  concerning total payoff criterion, are proved under much less restrictive assumptions than those used in the existing literature on single-population mean-field games. As single-popultaion games are just a specific case of the model presented here, in that way the paper also extends the theory for single-population mean-field games. This is further discussed when the relevant results are presented.

The organization of the paper is as follows: In Section 2 we present the way to model the discrete-time mean-field games with several populations of the players. In Section 3 we introduce some notation used in the remainder of the article. 
In Sections 4 and 5 we present several mean-field equilibrium-existence theorems for cases of discounted and total payoff, respectively. 
Finally, in Section 6 we give some concluding remarks.

\section{The model}
\label{sec:model}

Mean-field game models were designed to approximate dynamic game situations with a large number of symmetric agents. In multi-population mean-field games we still assume that the number of agents is large, but they are homogenic only within a smaller group called a {\em population}. The number of populations is finite and fixed, and their mutual interactions are encompassed in each individual's rewards and transitions. 
Each population has its own reward function and transition kernel (which may or may not operate on the same state space), which makes it significantly different from the models considered in the literature on discrete-time mean-field games. 
Below we describe the model formally. 

A {\em multi-population discrete-time mean-field game} is described by the following objects:
\begin{itemize}
\item We assume that the game is played in discrete time, that is
$t\in\{ 1,2,\ldots\}$.
\item The game is played by an infinite number (continuum) of
players divided into $N$ {\em populations}. Each player has a {\em private state} $s$, changing over
time. We assume that the set of individual states $S^i$ is the same for each player in population $i$ ($i=1,\ldots,N$), and that it is a nonempty closed subset of a locally compact Polish space $S$.\footnote{As it can be clearly seen, the model encompasses in particular the situation when the state space for each population is the same and equal to $S$.}
\item A vector $\overline{\mu}=(\mu^1,\ldots,\mu^N)\in\Pi_{i=1}^N\Delta(S^i)$ of $N$ probability distributions over Borel sets\footnote{Here and in the sequel, for any set $X$, $\Delta(X)$ denotes the set of probability distributions over the $\sigma$-algebra of Borel subsets of $X$, $\mathcal{B}(X)$.} of $S^i$, $i=1,\ldots,N$, is called a {\em global state} of the game. Its $i$-th component describes the proportion of $i$-th population, which is in each of the individual states.

We assume that at every stage of the game each player knows both his private state and the global state, and that his knowledge about individual states of his opponents is limited to the global state.
\item The set of {\em actions} available to a player from population $i$ in state $(s,\overline{\mu})$ is given by $A^i(s)$, 
with $A:=\bigcup_{i\in\{ 1,\ldots,N\}, s\in S^i}A^i(s)$ 
 -- a compact metric space. For any $i$, $A^i(\cdot)$ 
 is a non-empty compact valued
correspondence such that 
$$D^i:=\{ (s,a)\in S^i\times A: a\in A^i(s)\}$$
is a measurable set. Note that we assume that the sets of actions available to a player only depend on his private state and not on the global state of the game.
\item The global distribution of the state-action pairs is denoted by $\overline{\tau}=(\tau^1,\ldots,\tau^N)\in \Pi_{i=1}^N\Delta(D^i)$. Again, it gives the distributions of state-action pairs within the population divided into subpopulations $i=1,\ldots,N$.
\item {\em Immediate reward} of an individual from population $i$ is given by a measurable function $r^i:D^i\times\Pi_{i=1}^N\Delta(D^i)\rightarrow\mathbb{R}$. $r^i(s,a,\overline{\tau})$ gives the reward of a player at any stage of the game when his private state is $s$, his action is $a$ and the distribution of state-action pairs among the entire player population is $\overline{\tau}$.
\item {\em Transitions} are defined for each individual separately
with stochastic kernels $Q^i:D^i\times \Pi_{i=1}^N\Delta(D^i)\rightarrow \Delta(S^i)$ denoting transition probability for players from $i$-th population. $Q^i(B\mid\cdot,\cdot,\overline{\tau})$ is product-measurable for any $B\in\mathcal{B}(S^i)$, any $\overline{\tau}\in\Pi_{i=1}^N\Delta(D^i)$ and $i\in\{ 1,\ldots,N\}$.
\item The global state at time $t+1$, $\overline{\mu_t}$, is given by the aggregation of individual transitions of the players done by the formula
$$\mu_{t+1}^{i}(\cdot)=\Phi^i(\cdot\mid\overline{\tau_t}):=\int_{D^i}Q^i(\cdot\mid s,a,\overline{\tau_t})\tau_t^i(ds\times da).$$
As it can be clearly seen, the transition of the global state is deterministic.
\end{itemize}

A sequence $\pi^i=\{\pi_t^i\}_{t=0}^\infty$ of functions $\pi_t^i:S^i\rightarrow\Delta(A)$, such that $\pi_t^i(B\mid \cdot
)$ is measurable for any $B\in\mathcal{B}(A)$ and any $t$, 
satisfying $\pi_t^i(A^i(s)\mid s
)=1$ for every $s\in S^i$ and every $t$, is called a {\em Markov strategy} for a player of population $i$.
A function $f^i:S^i
\rightarrow \Delta(A)$, such that $f^i(B\mid \cdot
)$ is measurable for any $B\in\mathcal{B}(A)$, 
satisfying $f^i(A^i(s)\mid s
)=1$ for every $s\in S^i$ 
is called a {\em stationary strategy}. The set of all Markov strategies for players from $i$-th population is denoted by $\mathcal{M}^i$ while that of stationary strategies by $\mathcal{F}^i$. As in MDPs, stationary strategies can be seen as a specific case of Markov strategies that do not depend on $t$. In the paper we never consider general (history-dependent) strategies.\footnote{The lack of introduction of history-dependent strategies is only caused by the will to avoid additional complexity in the notation (which is rather complicated as it is). It should be noted though that, as in the case of discounted- and total-reward stochastic games, the Markov and stationary mean-field equilibria whose existence we prove later are in fact equilibria in the class of all strategies. This is a consequence of the fact that at a Markov (or stationary) equilibrium each agent faces the problem of maximizing a reward in a discounted (or total) reward Markov decision process which always admits a maximum for a Markov (stationary if the problem is time-homogeneous) policy.}

Next, let $\Pi^i_t(\pi^i,\mu^i)$ denote the state-action distribution of the $i$-th population players at time $t$ in the mean-field game corresponding to the distribution of individual states in population $i$, ${\mu}^i$ and a Markov strategy for players of population $i$, 
$\pi^i\in\mathcal{M}^i$, that is
$$\Pi_t^i(\pi^i,{\mu}^i)(B):=\int_B \pi_t^i(da\mid s)\mu^i(ds)\quad\mbox{for }B\in\mathcal{B}(D^i).$$
The vector $(\Pi_t^1(\pi^1,{\mu}^1),\ldots,\Pi_t^N(\pi^N,{\mu}^N))$ will be denoted by $\overline{\Pi}_t(\overline{\pi},\overline{\mu})$.
When we use this notation for stationary strategies, we skip the subscript $t$.

Given the evolution of the global state, which depends on the strategies of the players in a deterministic manner, we can define the individual history of a player $\alpha$ (from any given population $i$) as the sequence of his consecutive individual states and actions $h=(s^\alpha_0,a^\alpha_0,s^\alpha_1,a^\alpha_1,\ldots)$. By the Ionescu-Tulcea theorem (see Chap. 7 in \cite{BertsekasShreve}), for any Markov strategies $\pi^\alpha$ of player $\alpha$ and $\sigma^1,\ldots,\sigma^N$ of other players (including all other players of the same population), any initial global state $\overline{\mu_0}$ and any initial 
private state of player $\alpha$, $s$, 
there exists a unique probability measure $\mathbb{P}^{s,\overline{\mu_0},\overline{Q},\pi^\alpha,\overline{\sigma}}$ on the set of all infinite individual histories of the game $H=(D^i)^\infty$ endowed with Borel $\sigma$-algebra, such that for any $B\in\mathcal{B}(S^i)$, $E\in\mathcal{B}(A)$ and any partial history $h^\alpha_t=(s^\alpha_0,a^\alpha_0,\ldots,s^\alpha_{t-1},a^\alpha_{t-1},s^\alpha_t)\in (D^i)^t\times S^i=:H_t$, $t\in\mathbb{N}$,
\begin{equation}
\label{eq:1stIT}
\mathbb{P}^{s,\overline{\mu_0},\overline{Q},\pi^\alpha,\overline{\sigma}}(h\in H: s^\alpha_0=s)=1,
\end{equation}
\begin{equation}
\label{eq:2ndIT}
\mathbb{P}^{s,\overline{\mu_0},\overline{Q},\pi^\alpha,\overline{\sigma}}(h\in H: a^\alpha_t\in E\mid h^\alpha_t)=\pi_t^\alpha(E\mid s^\alpha_t),
\end{equation}
$$\mathbb{P}^{s,\overline{\mu_0},\overline{Q},\pi^\alpha,\overline{\sigma}}(h\in H: s^\alpha_{t+1}\in B\mid (h^\alpha_t,a^\alpha_t))=Q^i(B\mid s^\alpha_t,a^\alpha_t,\overline{\tau^t}),$$
with state-action distributions defined by $\tau^j_0=\Pi^j_0(\sigma^j,\mu^j_0)$, $\tau^j_{t+1}=\Pi^i_t(\sigma^j,\Phi^j(\cdot\mid \overline{\tau^t}))$ for $t=1,2,\ldots$ and $j=1,\ldots,N$.

Now we are ready to define the two types of reward we shall consider in this paper.
For $\beta\in(0,1)$, the {\em $\beta$-discounted reward}\footnote{Here we replace the superscript $\alpha$ used to define the measure $\mathbb{P}^{s,\overline{\mu_0},\overline{Q},\pi^\alpha,\overline{\sigma}}$ by $i$, as the situation is symmetric within the population.}  for a player $\alpha$ from population $i$ using
policy
$\pi^i\in\mathcal{M}^i$ 
when other players use policies
$\sigma^j\in\mathcal{M}^j$ (depending on the population $j$ they belong to) 
and the initial global state is $\overline{\mu_0}$, with the initial 
individual state of player $\alpha$ being $s^i_0$ is defined as follows:
$$J^i_\beta(s^i_0,\overline{\mu_0},\pi^i,\overline{\sigma})=\mathbb{E}^{s^i_0,\overline{\mu_0},\overline{Q},\pi^i,\overline{\sigma}}\sum_{t=0}^\infty\beta^tr^i(s_t^i,a_t^i,\overline{\tau^t}),$$
where  $\tau^j_0=\Pi_0^j(\sigma^j,\mu^j_0)$, $\tau^j_{t+1}=\Pi_t^j(\sigma^j,\Phi^j(\cdot\mid \overline{\tau^t}))$ for $t=1,2,\ldots$ and $j=1,\ldots,N$.


To define the total reward in our game let us distinguish one state
in $S$, say $s^*$, isolated from $S\setminus\{ s^*\}$ and assume that $A^i(s^*)=\{ a^*\}$
independently of $i\in\{ 1,\ldots,N\}$ for some fixed $a^*$ isolated from $A\setminus\{ a^*\}$. Moreover, let us assume that $s^*\in S^i$ for $i=1,\ldots,N$. Then the {\em total reward}
of a player from population $i$ using policy $\pi^i\in\mathcal{M}^i$ when other players
apply policies $\overline{\sigma}=(\sigma^1,\ldots,\sigma^N)$ and the initial global state is $\overline{\mu_0}$, with the initial 
individual state of player $\alpha$ being $s_0^i$, is defined in the
following way:
$$J^i_*(s^i_0,\overline{\mu_0},\pi^i,\overline{\sigma})=\mathbb{E}^{s^i_0,\overline{\mu_0},\overline{Q},\pi^i,\overline{\sigma}}\sum_{t=0}^{\mathcal{T}^i-1}r^i(s_t^i,a_t^i,\overline{\tau^t}),$$
where  $\tau^j_0=\Pi_0^j(\sigma^j,\mu^j_0)$, $\tau^j_{t+1}=\Pi_t^j(\sigma^j,\Phi^j(\cdot\mid \overline{\tau^t}))$ for $t=1,2,\ldots$ and $j=1,\ldots,N$,
while
$\mathcal{T}^i$ is the moment of the first arrival of the process
$\{ s_t^i\}$ to $s^*$. The total reward is interpreted as the reward accumulated by the
player over the whole of his lifetime. State $s^*$ is an artificial
state (so is action $a^*$) denoting that a player is dead. $\overline{\mu_0}$
corresponds to the distribution of the states across the population when he is
born, while $s^i_0$ is his own state when he is born. The fact that after some time the state of a
player can become again different from $s^*$ should be interpreted
as that after some time the player is replaced by some new-born one.


Next, we define the solutions we will be looking for:
\begin{definition}
Stationary strategies $f^1\in\mathcal{F}^1,\ldots f^N\in\mathcal{F}^N$ and a global state $\overline{\mu}\in\Pi_{i=1}^N\Delta(S^i)$ form a {\em stationary mean-field equilibrium} in the $\beta$-discounted reward game if for any $i$, $s^i_0\in S^i$, and every other stationary strategy of a player from population $i$,
$g^i\in\mathcal{F}^i$
$$J^i_\beta(s^i_0,\overline{\mu},f^i,\overline{f})\geq J^i_\beta(s^i_0,\overline{\mu},g^i,\overline{f})$$
and if $\overline{\mu}_0=\overline{\mu}$, then $\overline{\mu}_t=\overline{\mu}$ for every $t\geq 1$ if strategies $f^1,\ldots,f^N$ are used by all the players.

Markov strategies $\pi^1\in\mathcal{M}^1,\ldots \pi^N\in\mathcal{M}^N$ and a global state flow $(\overline{\mu}_0^*,\overline{\mu}_1^*,\ldots)\in(\Pi_{i=1}^N\Delta(S^i))^\infty$ form a {\em Markov mean-field equilibrium} in the $\beta$-discounted reward game if for any $i$, $s^i_0\in S^i$, and every other Markov strategy of a player from population $i$,
$\sigma^i\in\mathcal{M}^i$
$$J^i_\beta(s_0^i,\overline{\mu_0},\pi^i,\overline{\pi})\geq J^i_\beta(s_0^i,\overline{\mu_0},\sigma^i,\overline{\pi})$$
and if $\overline{\mu}_0=\overline{\mu}^*_0$ implies $\overline{\mu}_t=\overline{\mu}^*_t$ for every $t\geq 1$ if strategies $\pi^1,\ldots,\pi^N$ are used by all the players.
\end{definition}

Similarly,
\begin{definition}
Stationary strategies $f^1\in\mathcal{F}^1,\ldots f^N\in\mathcal{F}^N$ and a global state $\overline{\mu}\in\Pi_{i=1}^N\Delta(S^i)$ form a {\em stationary mean-field equilibrium} in the total reward game if for any $i$, $s_i^0\in S^i$, and every other stationary strategy of a player from population $i$,
$g^i\in\mathcal{F}^i$
$$J^i_*(s^i_0,\overline{\mu},f^i,\overline{f})\geq J^i_*(s^i_0,\overline{\mu},g^i,\overline{f}).$$
Moreover, if $\overline{\mu}_0=\overline{\mu}$, then $\overline{\mu}_t=\overline{\mu}$ for every $t\geq 1$ if strategies $f^1,\ldots,f^N$ are used by all the players.

Markov strategies $\pi^1\in\mathcal{M}^1,\ldots \pi^N\in\mathcal{M}^N$ and a global state flow $(\overline{\mu}_0^*,\overline{\mu}_1^*,\ldots)\in(\Pi_{i=1}^N\Delta(S^i))^\infty$ form a {\em Markov mean-field equilibrium} in the  total reward game if for any $i$, $t$, $s_i^t\in S^i$ and every other Markov strategy of a player from population $i$,
$\sigma^i\in\mathcal{M}^i$,
$$J^i_*(s^i_t,\overline{\mu}_t^*,{}^t\pi^i,{}^t\overline{\pi})\geq J^i_*(s^i_t,\overline{\mu}_t^*,{}^t\sigma^i,{}^t\overline{\pi}),$$
with ${}^ta$ denoting for any infinite vector $a=(a_0,a_1,\ldots)$, the vector $(a_t,a_{t+1},\ldots)$.
Moreover, if $\overline{\mu}_0=\overline{\mu}^*_0$, then $\overline{\mu}_t=\overline{\mu}^*_t$ for every $t\geq 1$ if strategies $\pi^1,\ldots,\pi^N$ are used by all the players.
\end{definition}

\section{Preliminaries}

As we have written, we assume that $S$ and $A$ are metric spaces. The metric on $S$ will be denoted by $d_S$ while that on $A$ by $d_A$. Whenever we relate to a metric on a product space, we mean the sum of the metrics on its coordinates.
Some of the assumptions presented below will be given with respect to the {\em moment function} $w_0:S\rightarrow [1,\infty)$, that is a continuous function satisfying
$$\lim_{n\rightarrow\infty}\inf_{s\in S\setminus K_n}w_0(s)=\infty$$
for some sequence $\{ K_n\} _{n\geq 1}$ of compact subsets of $S$. 
Moreover,
$$w_0(s)\geq 1+d_S(s,s_0)^p$$
for some $p\geq 1$ and $s_0\in S$.

In order to study both bounded and unbounded one-stage reward functions, we define the following function:
$$w:=\left\{ \begin{array}{ll} 1,&\mbox{ if each }r_i\mbox{ is bounded}\\w_0,&\mbox{ otherwise}\end{array}\right.$$
For any function $h:S\rightarrow\mathbb{R}$ we define its $w$-norm as
$$\left\| h\right\|_w:=\sup_{s\in S}\left\lvert\frac{h(s)}{w(s)}\right\rvert.$$
Whenever we speak of functions defined on a product of $S$ and some other space, their $w$-norm is defined similarly, with the help of the same function $w$.

By $B_w(S)$ we denote the space of all real-valued measurable functions from $S$ to $\mathbb{R}$ with finite $w$-norm. and by $C_w(S)$ -- the space of all continuous functions in $B_w(S)$. 
Clearly, both 
$B_w(S)$ and $C_w(S)$
are Banach spaces. The same can be said of $B_w(S\times A)$ and $C_w(S\times A)$ -- the spaces of bounded and bounded continuous functions  from $S\times A$ to $\mathbb{R}$ with finite $w$-norm.\footnote{We sometimes also use the notation $B_w(S^i)$, $C_w(S^i)$ for analogous sets of functions defined on these smaller domains.}

Analogously, for any finite signed measure $\mu$ on $S$, we define the $w$-norm of $\mu$ as
$$\left\|\mu\right\|_w=\sup_{g\in B_w(S),\| g\|_w\leq 1}\left\lvert \int_S g(s)\mu(ds)\right\rvert.$$
It should be noted that in case $w\equiv 1$, $\| \mu\|_w$ is the total variation distance (see e.g. \cite{HLL2}, Section 7.2.

There are two 
standard types of convergence of probability measures which are used in the paper: the weak convergence denoted by $\Rightarrow$ and the strong (or setwise) convergence denoted by $\rightarrow$ and defined (for any Borel space $(X,\mathcal{B}(X))$) by
$$\mu_n\rightarrow \mu\quad\Longleftrightarrow\quad \mu_n(B)\rightarrow\mu(B)\mbox{ for any }B\in\mathcal{B}(X).$$
It is known (see e.g. \cite{Parthasarathy}, Theorem 6.6) that the weak topology can be metrized using the metric
$$\rho(\mu,\nu):=\sum_{m=1}^\infty2^{-m}\left\lvert \int_S\phi_m(s)\mu(ds)-\int_S\phi_m(s)\nu(ds)\right\rvert,$$
where $\{\phi_i\}_{i\geq 1}$ is a sequence of continuous bounded functions from $S$ to $\mathbb{R}$ whose elements form a dense subset of the unit ball in $C(S)$. Strong convergence topology is in general not metric.

Next, let
$$\Delta_w(S):=\left\{ \mu\in\Delta(S):\int_Sw(s)\mu(ds)<\infty\right\}.$$
It has been shown in \cite{SaldiBasarRaginsky} that 
$\Delta_w(S)$ can be metrized using the metric
$$\rho_w(\mu,\nu):=\rho(\mu,\nu)+\left\lvert\int_Sw(s)\mu(ds)-\int_Sw(s)\nu(ds)\right\rvert$$
It can be shown that $\Delta_w(S)$ with metric $\rho_w$ is under the assumptions that we make about $w$ a Polish space (see \cite{SaldiBasarRaginsky,Bolley} for more on that). We will use the topology defined by this metric (called $w$-topology in the sequel) as the standard topology on $\Delta_w(S)$.

We will also use the notation
$$\Delta_w(S\times A):=\left\{ \tau\in\Delta(S\times A):\int_{S\times A}w(s)\tau(ds\times da)<\infty\right\}$$
with analogously defined metrics also denoted by $\rho$ (metric defining weak convergence) and $\rho_w$ ($w$-metric) as well as similar notation for subsets of $S$ or $S\times A$.

Whenever we speak about continuity of correspondences, we refer to the following definitions:\\
Let $X$ and $Y$ be two metric spaces and $F:X\rightarrow Y$, a correspondence. Let $F^{-1}(G)=\{ x\in X: F(x)\cap G\neq\emptyset\}$.
We say that $F$ is upper semicontinuous iff $F^{-1}(G)$ is closed for any closed $G\subset Y$. $F$ is lower semicontinuous iff $F^{-1}(G)$ is open for any open $G\subset Y$. $F$ is said to be continuous iff it is both upper and lower semicontinuous. For more on (semi)continuity of correspondences see \cite{HLL}, Appendix D or \cite{AliprantisBorder}, Chapter 17.2.

\section{The existence of stationary and Markov mean field equilibria in discounted payoff game}
\label{sec:mfg1}

\subsection{Assumptions}
\label{sec:ass}

In this section, we address the problem of the existence of an equilibrium in discrete-time mean-field games with $\beta$-discounted payoff. We begin by presenting the set of assumptions used in our results. 

\begin{enumerate}[{\bf ({A}1)}]
\item For $i=1,\ldots,N$, $r^i$ is continuous and bounded above by some constant $R$ on $D^i\times \Pi_{i=1}^N\Delta(D^i)$. Moreover,
for $i=1,\ldots,N$ and $s\in S^i$,
$$\inf_{(a,\overline{\tau})\in A^i(s)\times\Pi_{i=1}^N\Delta_w(D^i)}r^i(s,a,\overline{\tau})\geq -Rw(s).$$
\item 
For $i=1,\ldots,N$ and any sequence $\{ s_n,a_n,\overline{\tau}_n\}\subset D^i\times\Pi_{i=1}^N\Delta(D^i)$ such that $s_n\rightarrow s_*$, $a_n\rightarrow a_*$ and $\overline{\tau}_n\Rightarrow\overline{\tau}^*$, $Q^i(\cdot\mid s_n,a_n,\overline{\tau}_n)\rightarrow Q(\cdot\mid s_*,a_*,\overline{\tau}^*)$.
Moreover, 
\begin{enumerate}[(a)]
\item for $i=1,\ldots,N$ the functions
$$\int_S w(s')Q^i(ds'\mid s,a,\overline{\tau})$$ are continuous in $(s,a,\overline{\tau})$,
\item  for $i=1,\ldots,N$ and $s\in S^i$
$$\sup_{(a,\overline{\tau})\in A^i(s)\times\Pi_{i=1}^N\Delta(D^i)}\int_Sw(s')Q^i(s'\mid s,a,\overline{\tau})\leq w(s).$$
\end{enumerate}
\item For $i=1,\ldots,N$, correspondences $A^i$ are continuous.
\end{enumerate}

In some theorems weaker versions of assumptions (A1) and (A2) will be used:

\begin{enumerate}[{\bf ({A}1')}]
\item For $i=1,\ldots,N$, $r^i$ is continuous and bounded above by some constant $R$ on $D^i\times \Pi_{i=1}^N\Delta(D^i)$. Moreover,
there exist non-negative constants $\alpha$, $\gamma$, $M$ satisfying $\alpha\beta\gamma<1$ and 
$$\int_Sw(s)\mu_0^i(ds)\leq M\quad\mbox{for }i=1,\ldots,N,$$ 
and such that for $i=1,\ldots,N$, $s\in S^i$ and $t=0,1,2,\ldots$,
$$\inf_{(a,\overline{\tau})\in A^i(s)\times\Pi_{i=1}^N\Delta_w^{(t)}(D^i)}r^i(s,a,\overline{\tau})\geq -R\gamma^tw(s)$$
with $\Delta_w^{(t)}(D^i):=\left\{ \tau^i\in\Delta_w(D^i): \int_{D^i}w(s)\tau^i(ds\times da)\leq \alpha^tM\right\}$.
\item 
For $i=1,\ldots,N$ and any sequence $\{ s_n,a_n,\overline{\tau}_n\}\subset D^i\times\Pi_{i=1}^N\Delta_w(D^i)$ such that $s_n\rightarrow s_*$, $a_n\rightarrow a_*$ and $\overline{\tau}_n\Rightarrow\overline{\tau}^*$, $Q^i(\cdot\mid s_n,a_n,\overline{\tau}_n)\Rightarrow Q(\cdot\mid s_*,a_*,\overline{\tau}^*)$.
Moreover, 
\begin{enumerate}[(a)]
\item for $i=1,\ldots,N$ the functions
$$\int_S w(s')Q^i(ds'\mid s,a,\overline{\tau})$$ are continuous in $(s,a,\overline{\tau})$,
\item  for $i=1,\ldots,N$ and $s\in S^i$
$$\sup_{(a,\overline{\tau})\in A^i(s)\times\Pi_{i=1}^N\Delta(D^i)}\int_Sw(s')Q^i(s'\mid s,a,\overline{\tau})\leq\alpha w(s).$$
\end{enumerate}
\end{enumerate}

\subsection{Main results}

In the first two main results of this section we prove the existence of stationary mean-field equilibrium in 
discounted discrete-time mean-field games.

\begin{theorem}
\label{thm:strong_discounted_smfe_thm}
Suppose that the assumptions (A1--A3) are satisfied.
Then for any $\beta\in(0,1)$ the multi-population discrete-time mean-field game with $\beta$-discounted payoff defined with $r^i$, $Q^i$, $S^i$ and $A^i$, $i=1,\ldots,N$, has a stationary mean-field equilibrium.
\end{theorem}

In the proof we adapt the techniques introduced in \cite{JovanovicRosenthal} to our case.
We precede the proof of the theorem with two lemmas.
\begin{lemma}
\label{lem:V_cont}
For any $\overline{\tau}\in\Pi_{i=1}^N\Delta(D^i)$ let\footnote{The measure $\mathbb{P}^{s,\overline{Q},f^i}$ is defined here similarly as in the case of discounted rewards.}
$$V^i_{\beta,\overline{\tau}}(s):=\max_{f^i\in\mathcal{F}^i}\mathbb{E}^{\delta_s,\overline{Q},f^i}\sum_{t=0}^\infty\beta^tr^i(s_t^i,a_t^i,\overline{\tau}),$$
that is, let $V^i_{\beta,\overline{\tau}}$ be the optimal value for the $\beta$-discounted Markov decision process of player from population $i$ when the behaviour of all the other players is described by the state-action measure $\overline{\tau}$, fixed over time.
Under assumptions (A1--A3) $V^i_{\beta,\overline{\tau}}(s)$ is jointly continuous in $(s,\overline{\tau})$ on $S^i\times \Pi_{i=1}^N\Delta_w(D^i)$ and $\| V^i_{\beta,\overline{\tau}}\|_w\leq\frac{R}{1-\beta}$.
\end{lemma}

{\em Proof}:
Let us fix an $i\in\{ 1,\ldots,N\}$ and define for any $\overline{\tau}\in\Pi_{j=1}^N\Delta_w(D^j)$
$$T^i_{\overline{\tau}}(u)(s):=\sup_{a\in A^i(s)}\left[ r^i(s,a,\overline{\tau})+\beta\int_Su(s')Q^i(ds'\mid s,a,\overline{\tau})\right].$$
Note, that clearly (by assumptions (A1) and (A2) (b)) $T^i_{\overline{\tau}}$ maps $B_w(S^i)$ into itself. Moreover, for any $u_1,u_2\in B_w(S^i)$,
\begin{eqnarray}
\sup_{s\in S}\left\lvert \frac{T^i_{\overline{\tau}}(u_1)(s)-T^i_{\overline{\tau}}(u_2)(s)}{w(s)}\right\rvert&\leq&\sup_{s\in S, a\in A^i(s)}\frac{\beta\int_S\left\lvert(u_1(s')-u_2(s'))Q^i(ds'\mid s,a,\overline{\tau})\right\rvert}{w(s)}\nonumber\\
&\leq&\beta\sup_{s\in S,a\in A^i(s)}\frac{\int_S\|u_1-u_2\|_ww(s')Q^i(ds'\mid s,a,\overline{\tau})}{w(s)}\nonumber\\
&\leq&\beta\|u_1-u_2\|_w\sup_{s\in S
}\frac{w(s)}{w(s)}=\beta\|u_1-u_2\|_w,\label{eq:contraction}
\end{eqnarray}
where the penultimate inequality follows from the definition of the $w$-norm, while the last one from assumption (A2) (b).
Hence, $T^i_{\overline{\tau}}$ is a contraction defined on a complete space. By the Banach fixed point theorem it has a unique fixed point, which is by Theorem 4.2.3 in \cite{HLL} equal to $V^i_{\beta,\overline{\tau}}$. Moreover, this fixed point can be obtained as $\lim_{n\rightarrow\infty}\left( T^i_{\overline{\tau}}\right)^n(u_0)$ for any given $u_0\in B_w(S^i)$.

Let $u_0^{\overline{\tau}}\equiv 0$ and define for $n=1,2,\ldots$ $u_n^{\overline{\tau}}:=T^i_{\overline{\tau}}(u_{n-1}^{\overline{\tau}})$. We will next show that for each $n$, $u_n^{\overline{\tau}}(s)$ is continuous in $(s,{\overline{\tau}})$ and $\| u_n^{\overline{\tau}}\|_w\leq\frac{R}{1-\beta}$.

We prove these statements by induction on $n$.
For $n=0$ both claims are obvious. Suppose they hold for $n=k-1$.
Then by Theorem 3.3 in \cite{Serfozo} (see also Remark 3.4 (ii) there -- the assumptions given there are satisfied with $g=\frac{R}{1-\beta}w$ by (A2) (b)), $r^i(s,a,\overline{\tau})+\beta\int_Su_{k-1}^{\overline{\tau}}(s')Q^i(ds'\mid s,a,\overline{\tau})$ is jointly continuous in $(s,a,\overline{\tau})$, hence, by Proposition 7.32 in \cite{BertsekasShreve} 
$$u_{k}^{\overline{\tau}}(s)=\sup_{a\in A^i(s)}\left[ r^i(s,a,\overline{\tau})+\beta\int_Su_{k-1}^{\overline{\tau}}(s')Q^i(ds'\mid s,a,\overline{\tau})\right]$$
is also (jointly) continuous.
We also have (here the third inequality is a consequence of assumption (A2) (b), while the fourth one follows from (A1) and our inductive assumption):
\begin{align*}
&\left\| T^i_{\overline{\tau}}(u_{k-1}^{\overline{\tau}})\right\|_w\leq\sup_{(s,a)\in D^i}\frac{\lvert r^i(s,a,\overline{\tau})\rvert+\beta\left\lvert\int_Su_{k-1}^{\overline{\tau}}(s')Q^i(ds'\mid s,a,\overline{\tau})\right\rvert}{w(s)}\\
&\leq\sup_{(s,a)\in D^i}\frac{\lvert r^i(s,a,\overline{\tau})\rvert}{w(s)}+\sup_{(s,a)\in D^i}\frac{\beta\|u_{k-1}^{\overline{\tau}}\|_w\int_Sw(s')Q^i(ds'\mid s,a,\overline{\tau})}{w(s)}\\
&\leq\sup_{(s,a)\in D^i}\frac{\lvert r^i(s,a,\overline{\tau})\rvert}{w(s)}+\sup_{(s,a)\in D^i}\frac{\beta\|u_{k-1}^{\overline{\tau}}\|_ww(s)}{w(s)}\leq R+\frac{\beta R}{1-\beta}=\frac{R}{1-\beta}.
\end{align*}
Thus, the second claim has been proved for $n=k$.

To finish the proof, let us take convergent sequences $\{ s_k\}_{k\geq 1}$ in $S^i$ and $\{ \overline{\tau}_k\}_{k\geq 1}$ in $\Pi_{j=1}^N\Delta_w(D^j)$ such that $s_k\rightarrow s_*$ and $\overline{\tau}_k\Rightarrow \overline{\tau}_*$. We will show that
$V^i_{\beta,\overline{\tau}_k}(s_k)\rightarrow V^i_{\beta,\overline{\tau}_*}(s_*)$. 
We start the proof by noticing that the set $K:=\{ s_k: k\geq 1\}\cup\{ s_*\}$ is clearly compact, hence there exists a value $W$ such that $W\geq \lvert w(s)\rvert$ for $s\in K$.
Now, fix any $\varepsilon>0$ and let $n_0$ be such that
$$
W\beta^{n_0}\frac{R}{1-\beta}<\frac{\varepsilon}{3}.
$$

Clearly, by repeated use of (\ref{eq:contraction}) for $u_1=u_0^{\overline{\tau}_*}$ and $u_2=V^i_{\beta,\overline{\tau}_*}$, we obtain
$$\|u_{n_0}^{\overline{\tau}_*}-V^i_{\beta,\overline{\tau}_*}\|_w\leq\beta^{n_0}\left(\|u_{0}^{\overline{\tau}_*}-V^i_{\beta,\overline{\tau}_*}\|_w\right)=\beta^{n_0}\| V^i_{\beta,\overline{\tau}_*}\|_w\leq \beta^{n_0}\frac{R}{1-\beta},$$
hence
\begin{equation}
\label{eq:eps1}
\left\lvert u_{n_0}^{\overline{\tau}_*}(s_*)-V^i_{\beta,\overline{\tau}_*}(s_*)\right\rvert\leq w(s_*)\|u_{n_0}^{\overline{\tau}_*}-V^i_{\beta,\overline{\tau}_*}\|_w\leq W\beta^{n_0}\frac{R}{1-\beta}<\frac{\varepsilon}{3}.
\end{equation}
Similarly we obtain that for any $k\geq 1$,
\begin{equation}
\label{eq:eps2}
\left\lvert u_{n_0}^{\overline{\tau}_k}(s_k)-V^i_{\beta,\overline{\tau}_k}(s_k)\right\rvert\leq w(s_k)\|u_{n_0}^{\overline{\tau}_k}-V^i_{\beta,\overline{\tau}_k}\|_w\leq W\beta^{n_0}\frac{R}{1-\beta}<\frac{\varepsilon}{3}.
\end{equation}
Finally, from the joint continuity of $u_{n_0}^{\cdot}(\cdot)$, there exists a $k_0\in\mathbb{N}$ such that for any $k\geq k_0$
\begin{equation}
\label{eq:eps3}
\left\lvert u_{n_0}^{\overline{\tau}_k}(s_k)-u_{n_0}^{\overline{\tau}_*}(s_*)\right\rvert<\frac{\varepsilon}{3}.
\end{equation}

Now, combining (\ref{eq:eps1}), (\ref{eq:eps2}) and (\ref{eq:eps3}), we obtain that for any $k\geq k_0$
\begin{eqnarray*}
&&\left\lvert V^i_{\beta,\overline{\tau}_k}(s_k)-V^i_{\beta,\overline{\tau}_*}(s_*)\right\rvert\leq\left\lvert u_{n_0}^{\overline{\tau}_k}(s_k)-V^i_{\beta,\overline{\tau}_k}(s_k)\right\rvert\\
&&+\left\lvert u_{n_0}^{\overline{\tau}_k}(s_k)-u_{n_0}^{\overline{\tau}_*}(s_*)\right\rvert+\left\lvert u_{n_0}^{\overline{\tau}_*}(s_*)-V^i_{\beta,\overline{\tau}_*}(s_*)\right\rvert<\varepsilon,
\end{eqnarray*}
which ends the proof that $V^i_{\beta,\cdot}(\cdot)$ is continuous. Proof that $\|V^i_{\beta,\overline{\tau}}\|_w\leq\frac{R}{1-\beta}$ is elementary.
$\Box$

\begin{lemma}
\label{lem:inv_measure}
Suppose assumptions (A1--A3) hold and $M>0$ is such that for each $i\in \{ 1,\ldots,N\}$ the set 
$$\Delta_w^M(S^i):=\left\{ \mu\in\Delta(S^i): \int_{S^i}w(s)\mu(ds)\leq M\right\}$$ 
is nonempty. Then for each $\overline{\tau}\in\Pi_{j=1}^N\Delta(D^j)$, $i=1,\ldots,N$ and any stationary strategy $f^i\in\mathcal{F}^i$, there exists a $\mu_{f^i,\overline{\tau}}\in\Delta_w^M(S^i)$ such that
$$\mu_{f^i,\overline{\tau}}(B)
=\int_{S^i}\int_{A^i(s)}Q^i(B\mid s,a,\overline{\tau})f^i(da\mid s)\mu_{f^i,\overline{\tau}}(ds)$$
for any $B\in\mathcal{B}(S)$.
\end{lemma}

{\em Proof:}
Let us fix $i\in\{ 1,\ldots,N\}$ and note that for any $\overline{\tau}\in\Pi_{j=1}^N\Delta(D^j)$, and any stationary strategy $f^i\in\mathcal{F}^i$, the transition probability
\begin{equation}
\label{eq:Q_F_def}
Q^i(\cdot\mid s,f^i,\overline{\tau}):=\int_{A^i(s)}Q^i(\cdot\mid s,a,\overline{\tau})f^i(da\mid s)
\end{equation}
is clearly strongly continuous. 

Next, suppose that the initial distribution of individual state of a player from population $i$ is $\rho_0^i\in\Delta_w^M(S^i)$. We prove by induction that the same is true for $\rho_t^i$, the distribution of his state when $t=1,2,\ldots$ if he uses strategy $f^i$ and the behaviour of the other players is described by $\overline{\tau}$. Suppose the thesis is true for $t$. Then by assumption (A2) (b) we have
\begin{eqnarray*}
&&\int_{S^i}w(s)\rho_{t+1}^i(ds)=\int_{S^i}\int_{A^i(s)}\int_{S^i}w(s')Q^i(ds'\mid s,a,\overline{\tau})f^i(da\mid s)\rho_t^i(ds)\\
&\leq&\int_{S^i}\int_{A^i(s)}\left( \sup_{(a',\overline{\tau})\in A^i(s)\times\Pi_{i=1}^N\Delta(D^i)}\int_{S^i}w(s')Q^i(ds'\mid s,a',\overline{\tau})\right) f^i(da\mid s)\rho_t^i(ds)\\
&=&\int_{S^i}\left( \sup_{(a',\overline{\tau})\in A^i(s)\times\Pi_{i=1}^N\Delta(D^i)}\int_{S^i}w(s')Q^i(ds'\mid s,a',\overline{\tau})\right)\rho_t^i(ds)\\
&\leq&\int_{S^i}w(s)\rho_t^i(ds)\leq M
\end{eqnarray*}

By Remark 1 in \cite{HLL95} we know that the sequence of measures 
$$\nu^T(\cdot):=\frac{1}{T}\sum_{t=0}^{T-1}\rho_t^i(\cdot)$$
(whose elements clearly belong to $\Delta_w^M(S^i)$, as it is a convex set) has a subsequence weakly converging to an invariant measure of the Markov chain with transition probability $Q^i(\cdot\mid \cdot,f^i,\overline{\tau})$. Let us call it $\mu_{f^i,\overline{\tau}}$. It can be easily showed that $\Delta_w^M(S^i)$ is tight, hence, by Prohorov's theorem (Theorem 6.1 in \cite{Billingsley}) relatively compact. But this implies that $\Delta_w^M(S^i)$ is closed in weak convergence topology, hence $\mu_{f^i,\overline{\tau}}\in\Delta_w^M(S^i)$.
$\Box$

{\em Proof of Theorem \ref{thm:strong_discounted_smfe_thm}}: 
%
Let $M$ be such as in Lemma \ref{lem:inv_measure} and for $i=1,\ldots,N$ let
$$\Delta_w^M(D^i):=\left\{ \tau\in\Delta(D^i): \int_{D^i}w(s)\tau(ds\times da)\leq M\right\}.$$
Let us further define the correspondences from $\Pi_{j=1}^N\Delta_w^M(D^j)$ to $\Delta_w^M(D^i)$, $i=1,\ldots,N$:
$$
\Theta^i(\overline{\tau}):=\left\{ \eta^i\in\Delta_w^M(D^i): \eta^i_{S^i}(\cdot)=\int_{D^i}Q^i(\cdot\mid s,a,\overline{\tau})\eta^i(ds\times da)
\right\},
$$
\begin{eqnarray*}
\Psi^i_\beta(\overline{\tau})&:=&\left\{ \eta^i\in\Theta^i(\overline{\tau}): \int_{S^i}V^i_{\beta,\overline{\tau}}(s)\eta^i_{S^i}(ds)\right.\\&=&\left.\int_{D^i}\left[ r^i(s,a,\overline{\tau})+\beta\int_{S^i}V^i_{\beta,\overline{\tau}}(s')Q^i(ds'\mid s,a,\overline{\tau})\right]\eta^i(ds\times da)
\right\}
\end{eqnarray*}
We will now verify, that for each $i$, $\Theta^i$ and $\Psi^i$ have some useful properties. We fix $i\in\{ 1,\ldots,N\}$ for all these considerations.

First note that
$\eta=\Pi^i(f^i,\mu_{f^i,\overline{\tau}})$ clearly belongs to $\Theta^i(\overline{\tau})$, as for any $B\in\mathcal{B}(S^i)$,
\begin{eqnarray*}
&&\left( \Pi^i(f^i,\mu_{f^i,\overline{\tau}})\right)_{S^i}(B)=\mu_{f^i,\overline{\tau}}(B)
=\int_{S^i} Q^i(B\mid s,f^i,\overline{\tau})\mu_{f^i,\overline{\tau}}(ds)\\
&=&\int_{D^i} Q^i(B\mid s,a,\overline{\tau})f^i(da\mid s)\mu_{f^i,\overline{\tau}}(ds)=\int_{D^i}Q^i(B\mid s,a,\overline{\tau})\Pi^i(f^i,\mu_{f^i,\overline{\tau}})(ds\times da)
\end{eqnarray*}
where the first and the last equality  follow from the definition of $\Pi^i(\cdot,\cdot)$, the second comes from the definition of invariant measure, while the third one from (\ref{eq:Q_F_def}). Moreover, by Lemma \ref{lem:inv_measure} $\mu_{f^i,\overline{\tau}}\in\Delta_w^M(S^i)$, which immediately implies that $\eta=\Pi^i(f^i,\mu_{f^i,\overline{\tau}})\in\Delta_w^M(D^i)$.

We next show that 
the graph of $\Theta^i$ is closed in weak convergence topology. To prove that, first note that for any bounded continuous function $u:S^i\rightarrow\mathbb{R}$, $\int_{S^i}u(s)Q^i(ds\mid \cdot,\cdot,\cdot)$ is, by the strong continuity of $Q^i$, a continuous function, so for any sequences $\eta_n^i\in\Delta(D^i)$ and $\overline{\tau_n}\in\Pi_{j=1}^N\Delta(D^j)$ such that $\eta_n^i\in\Theta^i(\overline{\tau_n})$ with $\eta_n^i\Rightarrow\eta^i$ and $\overline{\tau_n}\Rightarrow\overline{\tau}$, $\int_{S^i}u(s)Q^i(ds\mid \cdot,\cdot,\overline{\tau_n})$ converges continuously to $\int_{S^i}u(s)Q^i(ds\mid \cdot,\cdot,\overline{\tau})$. Hence, by Theorem 3.3 in \cite{Serfozo}, we have
$$\int_{D^i}\int_{S^i}u(s)Q^i(ds\mid \widehat{s},\widehat{a},\overline{\tau_n})\eta_n^i(d\widehat{s}\times d\widehat{a})\rightarrow_{n\rightarrow\infty}\int_{D^i}\int_{S^i}u(s)Q^i(ds\mid \widehat{s},\widehat{a},\overline{\tau})\eta^i(d\widehat{s}\times d\widehat{a}),$$
which means that $\int_{D^i}Q^i(\cdot\mid s,a,\overline{\tau_n})\eta_n^i(ds\times da)\Rightarrow\int_{D^i}Q^i(\cdot\mid s,a,\overline{\tau})\eta^i(ds\times da)$. From the uniqueness of the limit this implies that $\left(\eta^i\right)_{S^i}=\int_{D^i}Q^i(\cdot\mid s,a,\overline{\tau})\eta^i(ds\times da)$, hence $\eta^i\in\Theta^i(\overline{\tau})$, which implies that the graph of $\Theta^i$ is closed.

By Theorem 4.2.3 in \cite{HLL} there exists an optimal stationary policy $f^i_*$ in the optimization problem of a player from population $i$ maximizing his discounted reward when the behaviour of all the other players is described by the global state $\overline{\tau}$, fixed over time. Moreover, $f^i_*$ is a measurable selector attaining maximum on the RHS of the equation
\begin{equation}
\label{eq:disc_optimality_eq}
V^i_{\beta,\overline{\tau}}(s)=\max_{a\in A^i(s)}\left[ r^i(s,a,\overline{\tau})+\beta\int_SV^i_{\beta,\overline{\tau}}(s')Q^i(ds'\mid s,a,\overline{\tau})\right].
\end{equation}
Then we can write that
\begin{eqnarray*}
&&\int_{S^i}V^i_{\beta,\overline{\tau}}(s)\left( \Pi^i(f^i_*,\mu_{f^i_*,\overline{\tau}})\right)_{S^i}(ds)=\int_{S^i}V^i_{\beta,\overline{\tau}}(s)\mu_{f^i_*,\overline{\tau}}(ds)\\
&=&\int_{S^i}\left[ r^i(s,f^i_*(s),\overline{\tau})+\beta\int_{S^i}V^i_{\beta,\overline{\tau}}(s')Q^i(ds'\mid s,f^i_*(s),\overline{\tau})\right] \mu_{f^i_*,\overline{\tau}}(ds)\\
&=&\int_{D^i}\left[ r^i(s,a,\overline{\tau})+\beta\int_{S^i}V^i_{\beta,\overline{\tau}}(s')Q^i(ds'\mid s,a,\overline{\tau})\right]\Pi^i(f^i_*,\mu_{f^i_*,\overline{\tau}})(ds\times da),
\end{eqnarray*}
which implies that $\Pi^i(f^i_*,\mu_{f^i_*,\overline{\tau}})\in\Psi^i_\beta(\overline{\tau})$.

Next we show that the graph of $\Psi^i_\beta$ is closed in the weak convergence topology. Let us take sequences $\overline{\tau_n}\in\Pi_{j=1}^N\Delta_w^M(D^j)$ and $\eta_n\in\Delta_w^M(D^i)$ such that $\eta_n\in\Theta(\overline{\tau_n})$ for every $n\in\mathbb{N}$ with $\eta_n\Rightarrow\eta$ and $\overline{\tau_n}\Rightarrow\overline{\tau}$. Since the graph of $\Theta^i$ is closed, to show that so is that of $\Psi^i_\beta$, we only need to prove that the equality defining the set $\Psi^i_\beta(\overline{\tau})$ is satisfied for $\overline{\tau}$ and $\eta$. Note however that for each $n$
\begin{eqnarray}
\label{eq:Bellman_n}
&&\int_{S^i}V^i_{\beta,\overline{\tau_n}}(s)(\eta_n)_{S^i}(ds)\\
&&=\int_{D^i}\left[ r^i(s,a,\overline{\tau_n})+\beta\int_{S^i}V^i_{\beta,\overline{\tau_n}}(s')Q^i(ds'\mid s,a,\overline{\tau_n})\right]\eta_n(ds\times da)\nonumber
\end{eqnarray}
Then by the continuity of $Q^i$ and $V^i_{\beta,\overline{\tau_n}}$ and Theorem 3.3 in \cite{Serfozo} (see also Remark 3.4 (ii) there -- by (A2) (b) the assumption presented there is true for $g=\frac{R}{1-\beta}w$), $\int_{S^i}V^i_{\beta,\overline{\tau_n}}(s')Q^i(ds'\mid \cdot,\cdot,\overline{\tau_n})$ converges continuously to $\int_{S^i}V^i_{\beta,\overline{\tau}}(s')Q^i(ds'\mid \cdot,\cdot,\overline{\tau})$. As also $r^i(\cdot,\cdot,\overline{\tau_n})$ converges continuously to $r^i(\cdot,\cdot,\overline{\tau})$ by (A1), we may apply Theorem 3.3 in \cite{Serfozo} again (again with $g=\frac{R}{1-\beta}w$, which satisfies the assumption given in Remark 3.4 (ii) by (A1) and (A2) (b)), we can pass to the limit in (\ref{eq:Bellman_n}), obtaining
$$\int_{S^i}V^i_{\beta,\overline{\tau}}(s)(\eta)_{S^i}(ds)=\int_{D^i}\left[ r^i(s,a,\overline{\tau})+\beta\int_{S^i}V^i_{\beta,\overline{\tau}}(s')Q^i(ds'\mid s,a,\overline{\tau})\right]\eta(ds\times da)$$
which ends the proof that the graph of $\Psi^i$ is closed.

Finally, we can also note that for each $\overline{\tau}$, the set $\Psi^i_\beta(\overline{\tau})$ is clearly convex.

Next, let us define the following correspondence mapping $\Pi_{i=1}^N\Delta(D^i)$ into itself:
$$\overline{\Psi}_\beta(\overline{\tau}):=\Psi^1_\beta(\overline{\tau})\times\ldots\times \Psi^N_\beta(\overline{\tau}).$$
It is obvious that our previous considerations imply that $\overline{\Psi}_\beta$ also has nonempty and convex values and that its graph is closed. To finish the proof we need to note that the function $w_0$ (which is a moment on $S$) is also a moment on $S\times A$ (as $A$ is compact), hence each $\Delta_w^M(D^i)$ is tight. Now Prohorov's theorem implies that $\Delta_w^M(D^i)$ is compact in weak convergence topology for $i=1,\ldots,N$ and $\Pi_{j=1}^N\Delta_w^M(D^j)$ is compact in product topology. Therefore,
by the Glickberg fixed point theorem \cite{Glicksberg52}, $\overline{\Psi}_\beta$ it has a fixed point.

Suppose $\overline{\tau}_*$ is this fixed point. 
By the well-known result, see e.g. \cite{Hinderer} p.89, for each $i\in\{ 1,\ldots,N\}$, $\tau^i_*$ 
can be disintegrated into a stochastic kernel $g^i_*\in\mathcal{F}^i_0$ and its marginal on $S^i$, $(\tau^i_*)_{S^i}$, that is, satisfying for any $D\in\mathcal{B}(D^i)$
\begin{equation}
\label{eq:disintegrate}
\tau^i_*(D)=\int_D g^i_*(da\mid s)(\tau^i_*)_{S^i}(ds).
\end{equation}
Let us further define 
$$S^i_0:=\left\{ s\in S^i :\int_{A^i(s)}\left[ r^i(s,a,\overline{\tau_*})+\beta\int_{S^i} V^i_{\beta,\overline{\tau_*}}(s')Q^i(ds'\mid s,a,\overline{\tau_*})\right]g^i_*(da\mid s)<V^i_{\beta,\overline{\tau_*}}(s)\right\} .$$
Then, since $\tau^i_*\in\Psi^i_\beta(\overline{\tau}_*)$, 
\begin{equation}
\label{eq:zero}
\tau^i_*(S^i_0)=0,
\end{equation} 
otherwise inequality in the definition of $S^i_0$ would imply an inequality in the definition of $\Psi^i_\beta$.

Let us thus define the strategy
$$\widehat{f^i}(s)=\left\{ \begin{array}{ll} g^i_*(s),&\mbox{ if }s\in S^i\setminus S^i_0\\
f^i_*(s),&\mbox{ if }s\in S^i_0\end{array}\right.$$
It is clear that for any $s\in S^i$,
$$V^i_{\beta,\overline{\tau}_*}(s)=\int_{A^i(s)}\left[ r^i(s,a,\overline{\tau}_*)+\beta\int_{S^i}V^i_{\beta,\overline{\tau}_*}(s')Q^i(ds'\mid s,a,\overline{\tau}_*)\right]\widehat{f^i}(da\mid s).$$
Then, for any $D\in\mathcal{B}(D^i)$ we can reason as follows:
\begin{eqnarray*}
\tau^i_*(D)&=&\int_D g^i_*(da\mid s)(\tau^i_*)_{S^i}(ds)=\int_{S^i}\int_{A^i(s)} \one_D(s,a)g^i_*(da\mid s)(\tau^i_*)_{S^i}(ds)\\
&=&\int_{S^i_0}\int_{A^i(s)} \one_D(s,a)g^i_*(da\mid s)(\tau^i_*)_{S^i}(ds)\\
&+&\int_{S^i\setminus S^i_0}\int_{A^i(s)} \one_D(s,a)g^i_*(da\mid s)(\tau^i_*)_{S^i}(ds)\\
&=&0+\int_{S^i\setminus S^i_0}\int_{A^i(s)} \one_D(s,a)\widehat{f^i}(da\mid s)(\tau^i_*)_{S^i}(ds)\\
&=&\int_{S^i_0}\int_{A^i(s)} \one_D(s,a)\widehat{f^i}(da\mid s)(\tau^i_*)_{S^i}(ds)\\
&+&\int_{S^i\setminus S^i_0}\int_{A^i(s)} \one_D(s,a)\widehat{f^i}(da\mid s)(\tau^i_*)_{S^i}(ds)\\
&=&\int_D\widehat{f^i}(da\mid s)(\tau^i_*)_{S^i}(ds),
\end{eqnarray*}
where the first equality follows from (\ref{eq:disintegrate}), the second and the last one -- from the definition of integral over a set, while the third and the fourth one use the definition of strategy $\widehat{f^i}$ and (\ref{eq:zero}).

Note however that the two last equalities proved imply that there exist strategies $\widehat{f^i}$ and invariant measures $\mu^i_*=(\tau^i_*)_S$ for each population $i\in\{ 1,\ldots,N\}$, such that $\widehat{f^i}$ is the best response in the $\beta$-discounted game against $\overline{\tau}_*$ and $\overline{\tau}_*$ is the stationary global state corresponding to the profile of strategies $(\overline{f^1},\ldots,\overline{f^N})$ and the initial global state $(\mu^1_*,\ldots,\mu^N_*)$, hence a stationary mean-field equilibrium in the $\beta$-discounted game.
$\Box$

\begin{remark}
It should be noted here that the results given in Theorem \ref{thm:strong_discounted_smfe_thm} applied to the model with a single population extend the existing results for such a case. The most general result of this type in the literature appears in \cite{JovanovicRosenthal} and concerns the case with compact individual state space.
\end{remark}

The last result in this section gives the conditions under which Markov mean-field equilibria in our models exist. They are based on one of the theorems given in \cite{SaldiBasarRaginsky}. It should be noted here that the assumptions in that paper are slightly stronger than in our model when applied to a single population. Namely, in our model the rewards and the transitions depend on the state-action distribution of the other players, while in \cite{SaldiBasarRaginsky} the dependence is only on the distribution of private states. Also, in our model we allow the set of feasible actions to depend on player's private state, while in \cite{SaldiBasarRaginsky} there was no such dependence. 
%

\begin{theorem}
\label{thm:strong_discounted_mmfe_thm}
Suppose that the assumptions (A1'--A2') and (A3) are satisfied.
Then for any $\beta\in(0,1)$ and any $\overline{\mu}_0\in\Pi_{j=1}^N\Delta_w(S^j)$ the multi-population discrete-time mean-field game with $\beta$-discounted payoff defined with $r^i$, $Q^i$, $S^i$ and $A^i$, $i=1,\ldots,N$, has a Markov mean-field equilibrium.
\end{theorem}

\begin{remark}
As we have already noted, in our model the rewards and the transitions may depend on the state-action distribution of the players, which differs from \cite{SaldiBasarRaginsky}, where the dependence is only on the distribution of private states. Such an assumption is not new to the mean-field game literature. In case of discrete-time games it has already been used in the first paper on this type of games \cite{JovanovicRosenthal}. It has also been applied in \cite{BB92,BB95,PiWEA,PiW_new}. As for the continuous time case, this type of models have been introduced by Gomes and Voskayan in \cite{GomesVoskayan} under the name of extended mean-field games. Cardaliaguet and Lehalle have proposed
the name of mean field games of controls in \cite{CardaliaguetLehalle} for this type of framework. Some further results on the topic include \cite{Kobeissi,AchdouKobeissi,CarmonaLacker,GPV,BHP,LauriereTangpi}.
\end{remark}

We precede the proof of Theorem \ref{thm:strong_discounted_mmfe_thm} by a counterpart of Lemma \ref{lem:V_cont} for the Markov case. It requires some additional notation. First, let us define for $i=1,\ldots,N$ the sets
$$\Xi^i:=\Pi_{t=0}^\infty \Delta_w^{(t)}(D^i),$$
$$\Xi:=\Pi_{i=1}^N\Xi^i.$$
Next, let for $t\geq 0$
$$L_t:=\sum_{k=t}^\infty(\alpha\beta)^{k-t}\gamma^kR=\frac{\gamma^tR}{1-\alpha\beta\gamma}.$$
Using these constants, we define for $i=1,\ldots,N$ and $t\geq 0$ the sets
$$\mathcal{C}_i^t:=\left\{ u\in C_w(S^i): \| u\|_w\leq L_t\right\},$$
$$\mathcal{C}_i:=\Pi_{t=0}^\infty\mathcal{C}_i^t.$$
It is easy to see that under (A1'), $\mathcal{C}_i$ with metric
$$\rho_{\mathcal{C}}((u_0,u_1,\ldots),(v_0,v_1,\ldots)):=\sum_{t=0}^\infty\delta^{-t}\| u_t-v_t\|_w,$$
where $\delta$ is chosen such that $\delta>\gamma$ and $\alpha\beta\delta<1$, is a complete metric space.

\begin{lemma}
\label{lem:V_cont_Markov}
For any state-action measure flow $\left(\overline{\tau}\right):=\left( \overline{\tau}_0,\overline{\tau}_1,\ldots\right)\in\Xi$, let
$$V^{i,t}_{\beta,\left(\overline{\tau}\right)}(s):=\max_{\pi^i\in\mathcal{M}^i}\mathbb{E}^{\delta_s,\overline{Q},\pi^i}\sum_{k=t}^\infty\beta^tr^i(s_k^i,a_k^i,\overline{\tau}_k),$$
that is, let it be the optimal value at time $t$ for the $\beta$-discounted Markov decision process of player from population $i$ when the behaviour of all the other players is described by the flow $\left(\overline{\tau}\right)$.
Under assumptions (A1'--A2') and (A3) for any $i\in\{ 1,\ldots,N\}$ and $t\geq 0$, $V^{i,t}_{\beta,\left(\overline{\tau}\right)}\in\mathcal{C}_i^t$.
\end{lemma}

{\em Proof}:
Let us fix an $i\in\{ 1,\ldots,N\}$ and define for any $\left(\overline{\tau}\right)\in\Xi$ and $t\geq 0$
$$T^{i,t}_{\left(\overline{\tau}\right)}(u)(s):=\sup_{a\in A^i(s)}\left[ r^i(s,a,\overline{\tau}_t)+\beta\int_Su(s')Q^i(ds'\mid s,a,\overline{\tau}_t)\right].$$
By Proposition 7.32 in \cite{BertsekasShreve} for any $u\in\mathcal{C}^i_{t+1}$, $T^{i,t}_{\left(\overline{\tau}\right)}(u)$ is continuous. Moreover,
\begin{eqnarray*}
\left\| T^{i,t}_{\left(\overline{\tau}\right)}(u)\right\|_w&\leq&\sup_{(s,a)\in D^i}\frac{\left\lvert r^i(s,a,\overline{\tau}_t)+\beta\int_Su(s')Q^i(ds'\mid s,a,\overline{\tau}_t)\right\rvert}{w(s)}\\
&\leq&\sup_{(s,a)\in D^i}\frac{R\gamma^tw(s)+\beta\alpha L_{t+1}w(s)}{w(s)}=R\gamma^t+\frac{\beta\alpha R\gamma^{t+1}}{1-\alpha\beta\gamma}=L_t,
\end{eqnarray*}
where the last inequality follows from (A1') and (A2') (b) (note that $\tau_t^i\in\Delta^{(t)}_w(D^i)$ is implied by the assumption that $\mu_0^i\in\Delta_w(D^i)$ and (b) of (A2') applied to the recursive formula for $\overline{\tau}_t$). Hence, $T^{i,t}_{\left(\overline{\tau}\right)}$ maps $\mathcal{C}^i_{t+1}$ into $\mathcal{C}^i_{t}$.
Next, for any $u_1,u_2\in \mathcal{C}^i_{t+1}$, we have
\begin{eqnarray}
\sup_{s\in S}\left\lvert \frac{T^{i,t}_{\left(\overline{\tau}\right)}(u_1)(s)-T^{i,t}_{\left(\overline{\tau}\right)}(u_2)(s)}{w(s)}\right\rvert&\leq&\sup_{s\in S, a\in A^i(s)}\frac{\beta\int_S\left\lvert(u_1(s')-u_2(s'))Q^i(ds'\mid s,a,\overline{\tau}_t)\right\rvert}{w(s)}\nonumber\\
&\leq&\beta\alpha\|u_1-u_2\|_w\sup_{s\in S
}\frac{w(s)}{w(s)}=\alpha\beta\|u_1-u_2\|_w,\label{eq:contraction2}
\end{eqnarray}
where the last inequality follows from the definition of the $w$-norm and the assumption (A2') (b).

We next define the operator $T^{i}_{\left(\overline{\tau}\right)}:\mathcal{C}_i\rightarrow\mathcal{C}_i$ with the formula
$$\left( T^{i}_{\left(\overline{\tau}\right)}(u_0,u_1,\ldots)\right)_t:=T^{i,t}_{\left(\overline{\tau}\right)}(u_{t+1})\quad\mbox{for }t\geq 0.$$
From what we have shown, it really maps $\mathcal{C}_i$ into itself. As (\ref{eq:contraction2}) implies that for any $(u_0,u_1,\ldots)$ and $(v_0,v_1,\ldots)$ in $\mathcal{C}_i$,
\begin{eqnarray*}
\rho_{\mathcal{C}}\left( T^{i}_{\left(\overline{\tau}\right)}(u_0,u_1,\ldots),T^{i}_{\left(\overline{\tau}\right)}(v_0,v_1,\ldots)\right)&\leq& \sum_{t=0}^\infty\delta^{-t}\alpha\beta\| u_{t+1}-v_{t+1}\|_w\\
&\leq&\alpha\beta\delta\rho_{\mathcal{C}}\left((u_0,u_1,\ldots),(v_0,v_1,\ldots)\right),
\end{eqnarray*}
it is an $\alpha\beta\delta$-contraction defined on a complete space. By the Banach fixed point theorem it has a unique fixed point. By Theorems 14.4 and 17.1 in \cite{Hinderer} the elements of this vector are equal to $V^{i,t}_{\beta,\left(\overline{\tau}\right)}$, $t\geq 0$, which ends the proof that the optimal value functions $V^{i,t}_{\beta,\left(\overline{\tau}\right)}\in\mathcal{C}_i^t$ for $t\geq 0$.
$\Box$

Now we are ready to pass to the main part of the proof of Theorem \ref{thm:strong_discounted_mmfe_thm}. 

{\em Proof of Theorem \ref{thm:strong_discounted_mmfe_thm}:}
We start by defining the correspondences from $\Xi$ into $\Xi^i$, ($i=1,\ldots,N$) with the formulas:
$$
\widetilde{\Theta}^i((\overline{\tau})):=\left\{ (\eta^i)\in\Xi^i\! : \left(\eta^i_0\right)_{S^i}\!=\!\mu_0^i\mbox{ and }\!\left(\eta^i_t\right)_{S^i}(\cdot)\!=\!\int_{D^i}\! Q^i(\cdot\mid s,a,\overline{\tau}_{t-1})\tau^i_{t-1}(ds\times da)
\right\}\!,
$$
\begin{align*}
&\widetilde{\Psi}^i_{\beta}((\overline{\tau})):=
\left\{ (\eta^i)\in\widetilde{\Theta}^i((\overline{\tau})): \int_{S^i}V^{i,t-1}_{\beta,(\overline{\tau})}(s)\left(\eta^i_t\right)_{S^i}(ds)\right.\\
&=\left.\int_{D^i}\left[ r^i(s,a,\overline{\tau}_{t-1})+\beta\int_{S^i}V^{i,t}_{\beta,(\overline{\tau})}(s')Q^i(ds'\mid s,a,\overline{\tau}_{t-1})\right]\eta^i_{t-1}(ds\times da)\mbox{ for }t\geq 1
\right\}.
\end{align*}
We next prove that $\widetilde{\Theta}^i$ and $\widetilde{\Psi}^i_{\beta}$ have some useful properties. We fix $i\in\{ 1,\ldots,N\}$ for these considerations. We start by showing that for any Markov strategy $\pi^i\in\mathcal{M}^i$ the flow $(\eta^i)$ defined with the recurrence
\begin{equation}
\label{eq:rec}
\eta^i_0:=\Pi_0^i(\pi^i,\mu_0^i),\quad\eta^i_{t+1}:=\Pi_t^i(\pi^i,\Phi^i(\cdot\mid \overline{\tau}_t))\mbox{ for }t=1,2,\ldots
\end{equation}
is an element of $\widetilde{\Theta}^i((\overline{\tau}))$. We do it by induction on $t$. For $t=0$ both $\left(\eta^i_0\right)_{S^i}=\mu_0^i$ and $\int_{D^i}w(s)\eta^i_0(ds\times da)\leq M$ are obvious (by the definition of $\Pi_0^i$ and Assumption (A1')). Now suppose 
$\int_{D^i}w(s)\eta^i_{t}(ds\times da)\leq M\alpha^{t}$. Then by the defnition of $\Phi^i$,
$$\left(\eta^i_{t+1}\right)_{S^i}(\cdot)=\Phi^i(\cdot\mid \overline{\tau}_t)=\int_{D^i}Q^i(\cdot\mid s,a,\overline{\tau}_{t})\tau^i_{t}(ds\times da).$$
Moreover, by (A2') (b) we have
\begin{eqnarray*}
&&\int_{D^i}w(s)\eta^i_{t+1}(ds\times da)=\int_{S^i}w(s)\Phi^i(ds\mid \overline{\tau}_t)\\
&=&\int_{S^i}\int_{A^i(s)}\int_{S^i}w(s')Q^i(ds'\mid s,a,\overline{\tau}_t)\tau_t^i(ds\times da)\\
&\leq&\int_{S^i}\int_{A^i(s)}\left( \sup_{(a',\overline{\sigma}_t)\in A^i(s)\times\Xi^i_t}\int_{S^i}w(s')Q^i(ds'\mid s,a',\overline{\sigma}_t)\right) \tau_t^i(ds\times da)\\
&=&\int_{S^i}\int_{S^i}\left( \sup_{(a',\overline{\sigma}_t)\in A^i(s)\times\Xi^i_t}\int_{S^i}w(s')Q^i(ds'\mid s,a',\overline{\sigma}_t)\right)\tau_t^i(ds\times da)\\
&\leq&\int_{S^i}\alpha w(s)\tau_t^i(ds\times da)\leq \alpha^{t+1}M,
\end{eqnarray*}
which shows that $\eta^i_{t+1}\in\Xi^i_{t+1}$ and by the induction principle that $(\eta^i)\in\Xi^i$.

Next we prove that the graph of $\widetilde{\Theta}^i$ is closed.
To do that, we take convergent sequences $\left\{(\overline{\tau}^{(n)})\right\}_{n\geq 0}$ in $\Xi$ and $\left\{(\eta^{i,(n)})\right\}_{n\geq 0}$ in $\Xi^i$ such that $(\eta^{i,(n)})\in\widetilde{\Theta}^i((\overline{\tau}^{(n)}))$ for each $n$. Moreover, $(\eta^{i,(n)})\Rightarrow(\eta^i)$ and $(\overline{\tau}^{(n)})\Rightarrow(\overline{\tau})$ as $n\rightarrow\infty$ for some $(\eta^i)\in\Xi^i$ and $(\overline{\tau})\in\Xi$. Now fix $t\geq 1$. By the joint strong continuity of $Q^i$
for any bounded continuous function $u:S^i\rightarrow\mathbb{R}$, $\int_{S^i}u(s)Q^i(ds\mid \cdot,\cdot,\overline{\tau}_{t-1}^{(n)})$ converges continuously to $\int_{S^i}u(s)Q^i(ds\mid \cdot,\cdot,\overline{\tau}_{t-1})$. Hence, by Theorem 3.3 in \cite{Serfozo}, we have
\begin{eqnarray*}
&&\int_{D^i}\int_{S^i}u(s)Q^i(ds\mid \widehat{s},\widehat{a},\overline{\tau}_{t-1}^{(n)})\tau^{i,(n)}_{t-1}(d\widehat{s}\times d\widehat{a})
\\
&&\rightarrow_{n\rightarrow\infty}\int_{D^i}\int_{S^i}u(s)Q^i(ds\mid \widehat{s},\widehat{a},\overline{\tau}_{t-1})\tau^i_{t-1}(d\widehat{s}\times d\widehat{a}),
\end{eqnarray*}
which means that $\int_{D^i}Q^i(\cdot\mid s,a,\overline{\tau}_{t-1}^{(n)})\tau^{i,(n)}_{t-1}(ds\times da)\Rightarrow\int_{D^i}Q^i(\cdot\mid s,a,\overline{\tau}_{t-1})\tau^i_{t-1}(ds\times da)$. Therefore, we have $\left(\eta^i_t\right)_{S^i}=\int_{D^i}Q^i(\cdot\mid s,a,\overline{\tau}_{t-1})\tau^i_{t-1}(ds\times da)$ for each $t\geq 0$, hence $\eta^i\in\widetilde{\Theta}^i((\overline{\tau}))$, which implies that the graph of $\widetilde{\Theta}^i$ is closed.


Next note, that if $\pi^i_\beta$ is an optimal deterministic Markov policy in the optimization problem of a player from population $i$ maximizing his $\beta$-discounted reward when the behaviour of all the other players at each stage is described by the flow $(\overline{\tau})$, then for each $t\geq 1$ and $s\in S^i$ it satisfies
$$V^{i,t-1}_{\beta,(\overline{\tau})}(s)=\left[ r^i(s,\pi^i_{\beta,t-1}(s),\overline{\tau}_{t-1})+\beta\int_{S^i}V^{i,t}_{\beta,(\overline{\tau})}(s')Q^i(ds'\mid s,\pi^i_{\beta,t-1}(s),\overline{\tau}_{t-1})\right],$$
which implies that for any $t\geq 1$,
\begin{align*}
&\int_{S^i}V^{i,t-1}_{\beta,(\overline{\tau})}(s)\left( \Pi^i_{t-1}(\pi^i_\beta,\Phi^i(\cdot\mid \overline{\tau}_{t-1}))\right)_{S^i}(ds)=\int_{S^i}V^{i,t-1}_{\beta,(\overline{\tau})}(s)\Phi^i(\cdot\mid \overline{\tau}_{t-1})(ds)\\
=&\int_{S^i}\left[ r^i(s,\pi^i_{\beta,t-1}(s),\overline{\tau}_{t-1})+\beta\int_{S^i}V^{i,t}_{\beta,(\overline{\tau})}(s')Q^i(ds'\mid s,\pi^i_{\beta,t-1}(s),\overline{\tau}_{t-1})\right] \Phi^i(\cdot\mid \overline{\tau}_{t-1})(ds)\\
=&\int_{D^i}\left[ r^i(s,a,\overline{\tau}_{t-1})+\beta\int_{S^i}V^{i,t}_{\beta,(\overline{\tau})}(s')Q^i(ds'\mid s,a,\overline{\tau}_{t-1})\right]\Pi^i_{t-1}(\pi^i_\beta,\Phi^i(\cdot\mid \overline{\tau}_{t-1}))(ds\times da),
\end{align*}
Hence, measure flow $(\eta^i)$ defined by (\ref{eq:rec}) with $\pi^i:=\pi^i_\beta$ is an element of $\widetilde{\Psi}^i_\beta((\overline{\tau}))$.

Next we prove that the graph of $\widetilde{\Psi}^i_\beta$ is closed.
Let $\left\{(\overline{\tau}^{(n)})\right\}_{n\geq 0}$ in $\Xi$ and $\left\{(\eta^{i,(n)})\right\}_{n\geq 0}$ in $\Xi^i$ be convergent sequences such  that $(\eta^{i,(n)})\in\widetilde{\Psi}^i_\beta((\overline{\tau}^{(n)}))$ for each $n$. Moreover, let $(\eta^{i,(n)})\Rightarrow(\eta^i)\in\Xi^i$ and $(\overline{\tau}^{(n)})\Rightarrow(\overline{\tau})\in\Xi$ as $n\rightarrow\infty$. 
Since the graph of $\widetilde{\Theta}^i$ is closed, showing that  $\widetilde{\Psi}^i_\beta$ has the same property only requires proving that the equalities defining $\widetilde{\Psi}^i_\beta$ hold for $(\eta^i)$ and $(\overline{\tau})$.
Let us fix $t\geq 1$.
The definition of $\widetilde{\Psi}^i_\beta$ implies that for each $n$
\begin{eqnarray}
\label{eq:betaBellman_n}
&&\int_{S^i}V^{i,t-1}_{\beta,(\overline{\tau}^{(n)})}(s)\left(\eta^{i,(n)}_t\right)_{S^i}(ds)\\
&=&\int_{D^i}\left[ r^i(s,a,\overline{\tau}^{(n)}_{t-1})+\beta\int_{S^i}V^{i,t}_{\beta,(\overline{\tau}^{(n)})}(s')Q^i(ds'\mid s,a,\overline{\tau}^{(n)}_{t-1})\right]\eta^{i,(n)}_{t-1}(ds\times da)\nonumber
\end{eqnarray}
By the continuity of $Q^i$ and $V^{i,t}_{\beta,(\overline{\tau}^{(n)})}$, and Theorem 3.3 in \cite{Serfozo} (the assumption presented in Remark 3.4 (ii) there is true for $g=L_tw$ by Lemma \ref{lem:V_cont_Markov}), $\int_{S^i}V^{i,t}_{\beta,(\overline{\tau}^{(n)})}(s')Q^i(ds'\mid \cdot,\cdot,\overline{\tau}^{(n)}_{t-1})$ converges continuously to $\int_{S^i}V^{i,t}_{\beta,(\overline{\tau}}(s')Q^i(ds'\mid \cdot,\cdot,\overline{\tau}_{t-1})$. As also $r^i(\cdot,\cdot,\overline{\tau}^{(n)}_{t-1})$ converges continuously to $r^i(\cdot,\cdot,\overline{\tau}_{t-1})$ by (A1'), using Theorem 3.3 in \cite{Serfozo} once more (now with $g=L_{t-1}w$ in Remark 3.4 (ii) there, again by by Lemma \ref{lem:V_cont_Markov}), we can pass to the limit in (\ref{eq:betaBellman_n}), obtaining
\begin{eqnarray*}
&&\int_{S^i}V^{i,t-1}_{\beta,(\overline{\tau})}(s)\left(\eta^i_t\right)_{S^i}(ds)\\
&=&\int_{D^i}\left[ r^i(s,a,\overline{\tau}_{t-1})+\beta\int_{S^i}V^{i,t}_{\beta,(\overline{\tau})}(s')Q^i(ds'\mid s,a,\overline{\tau}_{t-1})\right]\eta^i_{t-1}(ds\times da)
\end{eqnarray*}
As $t$ was arbitrary, this ends the proof that the graph of $\widetilde{\Psi}^i_\beta$ is closed.

To finalize the proof, we define the correspondence from $\Xi$ into itself:
$$\widetilde{\Psi}_\beta((\overline{\tau})):=\widetilde{\Psi}^1_\beta((\overline{\tau}))\times\ldots\times \widetilde{\Psi}^N_\beta((\overline{\tau})).$$
What we have shown already implies that $\widetilde{\Psi}_\beta$ has nonempty values and that its graph is closed. Convexity of values of $\widetilde{\Psi}_\beta$ is obvious. As $w$ is a moment function, each $\Delta_w^{(t)}(D^i)$ is tight, hence, by Prohorov's theorem it is compact. This implies that $\Xi$ is compact in product topology. Therefore,
by the Glickberg fixed point theorem \cite{Glicksberg52}, $\widetilde{\Psi}_\beta$ has a fixed point. Suppose $(\overline{\tau^*})$ is this fixed point. Disintegrating $\tau^{*i}_t$ gives for $i=1,\ldots,N$ and $t=0,1,\ldots$ stochastic kernels $\pi^{*i}_t$ and measures $\mu^{*i}_t$ which (after similar modifications as in the proof of Theorem \ref{thm:strong_discounted_smfe_thm}) correspond to the Markov strategies and global state flows in mean-field Markov equilibrium in our game.
$\Box$


\section{The existence of stationary and Markov mean field equilibria in total payoff game}
\label{sec:mfg2}

\subsection{Assumptions}
In this section, we address the problem of the existence of an equilibrium in the  mean-field games with total payoff. In its main results we shall add new assumption (A4) or (A4'') to those defined in Section \ref{sec:ass}.
Its formulation requires defining for $i=1,\ldots,N$, $s\in S^i$, $a\in A^i(s)$ and $\overline{\tau}\in\Pi_{j=1}^N\Delta(D^j)$ the modified transition probabilities $Q^i_*$:
$$Q^i_*(\cdot\mid s,a,\overline{\tau}):=\left\{ \begin{array}{ll}
Q^i(\cdot\mid s,a,\overline{\tau}),&\mbox{ if }s\neq s^*\\
\delta_{s^*},&\mbox{ if }s=s^*\end{array}\right.$$
The first new assumption will be used to prove the existence of a stationary mean-field equilibrium in total reward models.

\begin{enumerate}[{\bf ({A}4)}]
\item For $i=1,\ldots,N$ and $\overline{\tau}\in\Pi_{j=1}^N\Delta(D^j)$,
$$\lim_{T\rightarrow\infty}\sup_{\pi^i\in\mathcal{M}^i}\left\|\sum_{t=T}^\infty\int_{S^i\setminus\{ s^*\}}w(s')\left( Q^i_*\right)^t(ds'\mid s,\pi^i,\overline{\tau})\right\|_w=0.$$
\end{enumerate}
 
In case of the results about the existence of Markov mean-field equilibria in the discounted case, two of the assumptions (A1') and (A2') refered to the discount factor $\beta$, which does not exist in the total reward model. Hence, apart from the new assumption (A4''), new versions of these two assumptions will be necessary. For technical reasons, there are also some additional restrictions added.

\begin{enumerate}[{\bf ({A}1'')}]
\item For $i=1,\ldots,N$, $r^i$ is continuous and bounded above by some constant $R$ on $D^i\times \Pi_{i=1}^N\Delta(D^i)$. Moreover,
there exist non-negative constants $\alpha$, $\gamma$, $M$ satisfying $\alpha\leq \gamma$, $\alpha\gamma<1$ and 
$$\int_Sw(s)\mu_0^i(ds)\leq M\quad\mbox{for }i=1,\ldots,N,$$ 
and such that for $i=1,\ldots,N$, $s\in S^i$ and $t=0,1,2,\ldots$,
$$\inf_{(a,\overline{\tau})\in A^i(s)\times\Pi_{i=1}^N\Delta_w^{(t)}(D^i)}r^i(s,a,\overline{\tau})\geq -R\gamma^tw(s)$$
with $\Delta_w^{(t)}(D^i)=\left\{ \tau^i\in\Delta_w(D^i): \int_{D^i}w(s)\tau^i(ds\times da)\leq \alpha^tM\right\}$.
\item 
For $i=1,\ldots,N$ and any sequence $\{ s_n,a_n,\overline{\tau}_n\}\subset D^i\times\Pi_{i=1}^N\Delta_w(D^i)$ such that $s_n\rightarrow s_*$, $a_n\rightarrow a_*$ and $\overline{\tau}_n\Rightarrow\tau^*$, $Q^i(\cdot\mid s_n,a_n,\overline{\tau}_n)\rightarrow Q(\cdot\mid s_*,a_*,\tau^*)$.
Moreover, 
\begin{enumerate}[(a)]
\item for $i=1,\ldots,N$ the functions
$$\int_S w(s')Q^i(ds'\mid s,a,\overline{\tau})$$ are continuous in $(s,a,\overline{\tau})$,
\item  for $i=1,\ldots,N$ and $s\in S^i$
$$\sup_{(a,\overline{\tau})\in A^i(s)\times\Pi_{i=1}^N\Delta(D^i)}\int_Sw(s')Q^i(s'\mid s,a,\overline{\tau})\leq\alpha w(s).$$
\end{enumerate}
\end{enumerate}

\begin{enumerate}[{\bf ({A}4'')}]
\item For $i=1,\ldots,N$,
$$\lim_{T\rightarrow\infty}\sup_{\substack{\pi^i\in\mathcal{M}^i,\\(\overline{\tau})\in\Pi_{t=0}^\infty\Pi_{j=1}^N\Delta(D^j)}}\left\|\sum_{t=T}^\infty\int_{S^i\setminus\{ s^*\}}w(s')\alpha^{-t}\left( Q^i_*\right)^t(ds'\mid s,\pi^i,(\overline{\tau}))\right\|_w=0.$$
\end{enumerate}

Here and in the sequel for $i=1,\ldots,N$, $s\in S^i$, $\pi^i\in\mathcal{M}^i$ and $(\overline{\tau})=(\overline{\tau}_0,\overline{\tau}_1,\ldots)\in\Pi_{t=0}^\infty\Pi_{j=1}^N\Delta(D^j)$,
$$\left( Q^i_*\right)^1(\cdot\mid s,\pi^i,(\overline{\tau})):=\int_{a\in A^i(s)}Q^i_*(\cdot\mid s,a,\overline{\tau}_0)\pi_0^i(da\mid s)$$
and for $t\geq 2$:
$$\left( Q^i_*\right)^t(\cdot\mid s,\pi^i,(\overline{\tau})):=\int_{S^i}\int_{A^i(s)}Q^i(\cdot\mid s',a,\overline{\tau}_t)\pi_{t-1}^i(da\mid s')\left( Q^i_*\right)^{t-1}(ds'\mid s,\pi^i,(\overline{\tau})).$$
For $i=1,\ldots,N$, $s\in S^i$, $\pi^i\in\mathcal{M}^i$ and $\overline{\tau}\in\Pi_{j=1}^N\Delta(D^j)$,
$$\left( Q^i_*\right)^t(\cdot\mid s,\pi^i,\overline{\tau}):=\left( Q^i_*\right)^t(\cdot\mid s,\pi^i,(\overline{\tau},\overline{\tau},\ldots)).$$

\begin{remark}
By assuming (A4) or (A4''), we build upon the framework of transient total-reward Markov decision processes introduced by Veinott \cite{Veinott} in the context of finite state and action spaces and generalized to Borel spaces in \cite{Pliska,FeinbergHuang}. The optimization problem faced by an individual from population $i$ in our model of total reward mean-field game would for every fixed global state-action distribution $\overline{\tau}$ be transient, if
$$\sup_{\pi^i\in\mathcal{M}^i}\left\|\sum_{t=0}^\infty\int_{S^i\setminus\{ s^*\}}w(s')\left( Q^i_*\right)^t(ds'\mid s,\pi^i,\overline{\tau})\right\|_w<\infty,$$
which is clearly true under (A2) (b) and (A4).
Roughly speaking, this would mean that that for a reward function such that $\left\| r^i(\cdot,\cdot,\overline{\tau})\right\|_w<\infty$, the total reward is finite for any Markov strategy applied. In (A4) we strengthen this assumption by requireing that the convergence of the total reward to its value is uniform across all Markov strategies with respect to the $w$-norm. (A4'') is an adjustement of this condition to the case when the decision-maker optimizes his behavior against a flow $(\overline{\tau})$.
\end{remark}

\subsection{Main results}

\begin{theorem}
\label{thm:strong_total_smfe_thm}
Suppose that the assumptions (A1--A4) are satisfied.
Then the multi-population discrete-time mean-field game with total payoff defined with $r^i$, $Q^i$, $S^i$ and $A^i$, $i=1,\ldots,N$, has a stationary mean-field equilibrium.
\end{theorem}

Let us start by noticing that total reward of a player from population $i$ using a given strategy $\pi$ when the behaviour of the others is constant over time and described by $\overline{\tau}$ in the MDP model with transition probability $Q^*_i$ is the same as the reward until reaching state $s^*$ in the model with transition probability $Q_i$. Let us next define for any $i\in\{ 1,\ldots,N\}$ and $\overline{\tau}\in\Pi_{j=1}^N\Delta(D^j)$
$$V^i_{*\overline{\tau}}(s):=\max_{f^i\in\mathcal{F}^i}\mathbb{E}^{\delta_s,\overline{Q_*},f^i}\sum_{t=0}^\infty r^i(s_t^i,a_t^i,\overline{\tau}),$$
that is, the optimal value for the total-reward Markov decision process of player from population $i$ when the behaviour of all the other players is described by the state-action measure $\overline{\tau}$, fixed over time. Crucial properties of function $V^i_{*\cdot}(\cdot)$ are given in lemma below.

\begin{lemma}
\label{lem:total_V_cont}
Under assumptions (A1--A4) for each $i\in\{ 1,\ldots,N\}$, $V^i_{*\overline{\tau}}(s)$ is jointly continuous in $(s,\overline{\tau})$. Moreover, there exists a constant $L$ such that $\left\|V^i_{*\overline{\tau}}(\cdot)\right\|_w\leq L$ for any $\overline{\tau}\in\Pi_{j=1}^N\Delta_w(D^j)$.
\end{lemma}

{\em Proof}: Fix $i$ and $\overline{\tau}$, and note that under (A2) (b) and (A4),
$$L_i:=\left\|\sum_{t=0}^\infty\int_{S^i\setminus\{ s^*\}}w(s')\left( Q^i_*\right)^t(ds'\mid s,\pi^i,\overline{\tau})\right\|_w$$
is finite for $i=1,\ldots,N$. Hence, we immediately see that
the total-reward Markov decision process defined with $S^i$, $A^i$, $r^i$ and $Q^i_*$ satisfies the assumptions\footnote{With weight function $V=w$, $K=L_i$ and $\overline{c}=R$.} of Theorem 12 in \cite{FeinbergHuang}. Therefore, with a help of this theorem and Proposition 1 in \cite{FeinbergHuang} we can define:
\begin{enumerate}[(a)]
\item The function $\zeta^i:S^i\times\Pi_{j=1}^N\Delta_w(D^j)\rightarrow[0,\infty)$
\begin{equation}
\label{eq:mu_function}
\zeta^i(s,\overline{\tau}):=\sup_{\pi^i\in\mathcal{M}^i}\left[\sum_{t=0}^\infty\int_{S^i\setminus\{ s^*\}}w(s')\left( Q^i_*\right)^t(ds'\mid s,\pi^i,\overline{\tau})\right]
\end{equation}
satisfying for $s\in S^i$ and $\overline{\tau}\in\Pi_{j=1}^N\Delta_w(D^j)$, $\zeta^i(s,\overline{\tau})\in [w(s),L_iw(s)]$,
\item The discount factor $\beta:=\max_{j=1,\ldots,N}\frac{L_j-1}{L_j}$,
\item Modified one-stage rewards $r^i_*$:
$$r^i_*(s,a,\overline{\tau}):=\frac{r^i(s,a,\overline{\tau})}{\zeta^i(s,\overline{\tau})},$$
\item Modified transition probabilities $Q^i_{**}$ given by
\begin{align*}
&Q^i_{**}(B\mid s,a,\overline{\tau}):=\\
&\left\{ \begin{array}{ll}
\frac{1}{\beta\zeta^i(s,\overline{\tau})}\int_{B}\zeta^i(s',\overline{\tau})Q^i_*(ds'\mid s,a,\overline{\tau}),&\mbox{ if }B\in\mathcal{B}(S^i\setminus\{ s^*\}),\\
& (s,a)\in D^i\setminus\{ (s^*,a^*)\}\\
1-\frac{1}{\beta\zeta^i(s,\overline{\tau})}\int_{B}\zeta^i(s',\overline{\tau})Q^i_*(ds'\mid s,a,\overline{\tau}),&\mbox{ if }B=\{ s^*\}),\\ &(s,a)\in D^i\setminus\{ (s^*,a^*)\}\\
1,&\mbox{ if }B=\{ s^*\}, (s,a)=(s^*,a^*)
\end{array}\right.
\end{align*}
\end{enumerate}
such that $S^i$, $A^i$, $r^i_*$ and $Q^i_{**}$ define a $\beta$-discounted Markov decision process with a value 
\begin{equation}
\label{V_beta_vs_V_*}
V^i_{\beta\overline{\tau}}(s)=\frac{V^i_{*\overline{\tau}}(s)}{\zeta(s,\overline{\tau})}
\end{equation} 
for $s\in S^i$, $\overline{\tau}\in\Pi_{j=1}^N\Delta_w(D^j)$. Moreover, optimal stationary strategies exist and coincide in both MDPs.

Next note, that if $\zeta^i$ is a continuous function, then the model defined by $S^i$, $A^i$, $r^i_*$ and $Q^i_{**}$ satisfies assumptions (A1--A3) with function $w$ replaced by $w_*\equiv 1$ (in particular, (A2) follows from the fact that $\zeta^i(s,\cdot)\geq w(s)$ for $s\in S^i$). This implies that by Lemma \ref{lem:V_cont}, $V^i_{\beta\overline{\tau}}(s)$ is continuous in $(s,\overline{\tau})$ and $\left\|V^i_{\beta\overline{\tau}}(\cdot)\right\|\leq\frac{R}{1-\beta}$. Hence, combining (\ref{V_beta_vs_V_*}) with the fact that $\zeta^i$ is continuous and $\zeta^i(s,\cdot)\leq L_iw(s)$ we obtain that $V^i_{*\overline{\tau}}(s)$ is also continuous in $(s,\overline{\tau})$. Moreover, for any $\overline{\tau}\in\Pi_{j=1}^N\Delta_w(D^j)$, 
$$\left\| V^i_{*\overline{\tau}}(\cdot)\right\|_w\leq \frac{L_iR}{1-\beta}\leq \frac{\max_{j=1,\ldots,N}L_jR}{1-\beta}=:L,$$ 
which proves the thesis of the lemma. Hence, all we need to do is to show that $\zeta^i$ is continuous.

To do that, we note that $\zeta^i$ is clearly the limit of the sequence of functions $\left\{ w_n^{\overline{\tau}}\right\}_{n\geq 0}$, defined with the following recurrence: $w_0^{\overline{\tau}}:=w$, $w_n^{\overline{\tau}}:=T^i_{*\overline{\tau}}(w_{n-1}^{\overline{\tau}})$ for $n=1,2,\ldots$, where for any $\overline{\tau}\in\Pi_{j=1}^N\Delta_w(D^j)$,
$$T^i_{*\overline{\tau}}(u)(s):=\sup_{a\in A^i(s)}\left[ w(s)+\int_{S^i\setminus\{ s^*\}}u(s')Q^i_*(ds\mid s,a,\overline{\tau})\right].$$
We next show by indution that each $w_n^{\overline{\tau}}(s)$ is continuous in $(s,\overline{\tau})$. For $n=0$ the claim is true by the definition of $w$. Suppose it holds for $n=k-1$.
Then by Theorem 3.3 in \cite{Serfozo} (the assumptions given in Remark 3.4 (ii) there are satisfied with $g=L_iw$ because $w_{k-1}^{\overline{\tau}}(\cdot)\leq \zeta^i(\cdot,\overline{\tau})\leq L_iw(\cdot)$), $w(s)+\int_{S^i\setminus\{ s^*\}}w_{k-1}^{\overline{\tau}}(s')Q^i_*(ds'\mid s,a,\overline{\tau})$ is jointly continuous in $(s,a,\overline{\tau})$, hence, by Proposition 7.32 in \cite{BertsekasShreve} 
$$w_{k}^{\overline{\tau}}(s)=\sup_{a\in A^i(s)}\left[ w(s)+\int_{S^i\setminus\{ s^*\}}w_{k-1}^{\overline{\tau}}(s')Q^i_*(ds'\mid s,a,\overline{\tau})\right]$$
is also (jointly) continuous. Therefore, the claim is true for any $n\geq 1$.

To finish the proof, let us take convergent sequences $\{ s_k\}_{k\geq 1}$ in $S^i$ and $\{ \overline{\tau}_k\}_{k\geq 1}$ in $\Pi_{j=1}^N\Delta_w(D^j)$ such that $s_k\rightarrow s_*$ and $\overline{\tau}_k\Rightarrow \overline{\tau}_*$. We will show that
$\zeta^i(s_k,\overline{\tau}_k)\rightarrow \zeta^i(s_*,\overline{\tau}_*)$. 
Since the set $K:=\{ s_k: k\geq 1\}\cup\{ s_*\}$ is clearly compact, there exists a value $W$ such that $W\geq \lvert w(s)\rvert$ for $s\in K$.
Now, fix any $\varepsilon>0$. By (A4) there exists an $t^*$ such that
$$\sup_{\pi^i\in\mathcal{M}^i}\left\|\sum_{t=t^*+1}^\infty\int_{S^i\setminus\{ s^*\}}w(s')\left( Q^i_*\right)^t(ds'\mid s,\pi^i,\overline{\tau})\right\|_w<\frac{\varepsilon}{3W}.$$

This immediately implies
\begin{align}
&\left\lvert w_{t^*}^{\overline{\tau}_*}(s_*)-\zeta^i(s_*,\overline{\tau}_*)\right\rvert\leq
\left\lvert\sup_{\pi^i\in\mathcal{M}^i}\sum_{t=t^*+1}^\infty\int_{S^i\setminus\{ s^*\}}w(s')\left( Q^i_*\right)^t(ds'\mid s_*,\pi^i,\overline{\tau})\right\rvert\nonumber\\
\label{eq:eps1t}
&\leq
w(s_*)\sup_{\pi^i\in\mathcal{M}^i}\left\|\sum_{t=t^*+1}^\infty\int_{S^i\setminus\{ s^*\}}w(s')\left( Q^i_*\right)^t(ds'\mid s_*,\pi^i,\overline{\tau})\right\|_w
\leq W\frac{\varepsilon}{3W}=\frac{\varepsilon}{3}.
\end{align}
and, for any $k\geq 1$,
\begin{align}
&\left\lvert w_{t^*}^{\overline{\tau}_k}(s_k)-\zeta^i(s_k,\overline{\tau}_k)\right\rvert\nonumber\\
&\leq w(s_k)\sup_{\pi^i\in\mathcal{M}^i}\left\|\sum_{t=t^*+1}^\infty\int_{S^i\setminus\{ s^*\}}w(s')\left( Q^i_*\right)^t(ds'\mid s_k,\pi^i,\overline{\tau})\right\|_w\leq\frac{\varepsilon}{3}
\label{eq:eps2t}
\end{align}
Finally, from the joint continuity of $w_{t^*}^{\cdot}(\cdot)$, there exists a $k_0\in\mathbb{N}$ such that for any $k\geq k_0$
\begin{equation}
\label{eq:eps3t}
\left\lvert w_{t^*}^{\overline{\tau}_k}(s_k)-w_{t^*}^{\overline{\tau}_*}(s_*)\right\rvert<\frac{\varepsilon}{3}.
\end{equation}

Combining (\ref{eq:eps1t}), (\ref{eq:eps2t}) and (\ref{eq:eps3t}), we obtain that for any $k\geq k_0$
\begin{eqnarray*}
&&\left\lvert \zeta^i(s_k,\overline{\tau}_k)-\zeta^i(s_*,\overline{\tau}_*)\right\rvert\leq\left\lvert \zeta^i(s_k,\overline{\tau}_k)-w_{t^*}^{\overline{\tau}_k}(s_k)\right\rvert\\
&&+\left\lvert w_{t^*}^{\overline{\tau}_k}(s_k)-w_{t^*}^{\overline{\tau}_*}(s_*)\right\rvert+\left\lvert w_{t^*}^{\overline{\tau}_*}(s_*)-\zeta^i(s_*,\overline{\tau}_*)\right\rvert<\varepsilon,
\end{eqnarray*}
which ends the proof that $\zeta^i(\cdot,\cdot)$ is continuous.
$\Box$

{\em Proof of Theorem \ref{thm:strong_total_smfe_thm}}: 
As in the case of the discounted reward, we define the correspondences from $\Pi_{j=1}^N\Delta(D^j)$ to $\Delta(D^i)$:
$$
\Theta^i(\overline{\tau}):=\left\{ \eta^i\in\Delta_w^M(D^i): \eta^i_{S^i}(\cdot)=\int_{D^i}Q^i(\cdot\mid s,a,\overline{\tau})\eta^i(ds\times da)
\right\},
$$
\begin{eqnarray*}
\Psi^i_*(\overline{\tau})&:=&\left\{ \eta^i\in\Theta^i(\overline{\tau}): \int_{S^i}V^i_{*\overline{\tau}}(s)\eta^i_{S^i}(ds)\right.\\&=&\left.\int_{D^i}\left[ r^i(s,a,\overline{\tau})+\int_{S^i}V^i_{*\overline{\tau}}(s')Q^i_*(ds'\mid s,a,\overline{\tau})\right]\eta^i(ds\times da)
\right\}
\end{eqnarray*}
Using similar arguments to those employed in the proof of Theorem \ref{thm:strong_discounted_smfe_thm} (with a difference that instead of Lemma \ref{lem:V_cont} we apply Lemma \ref{lem:total_V_cont} when necessary), we can prove that $\Psi^i_*$ (we do not need to prove that for $\Theta^i$, as it is defined in exactly the same way as in the case of the $\beta$-discounted reward) that it has nonempty convex values and that its graph is closed. Then we define the correspondence
$$\overline{\Psi}_*(\overline{\tau}):=\Psi^1_*(\overline{\tau})\times\ldots\times \Psi^N_*(\overline{\tau})$$
and show that it has a fixed point, which, again using similar arguments as in the proof of Theorem \ref{thm:strong_discounted_smfe_thm}, can be proved to correspond to a stationary mean-field equilibrium in the total reward discrete-time mean field game considered in the theorem. 
$\Box$

In the last result of this section we give conditions under which a Markov mean-field equilibrium exists in the total-reward game.
\begin{theorem}
\label{thm:strong_total_mmfe_thm}
Suppose that the assumptions (A1''), (A2''), (A3) and (A4'') are satisfied.
Then for any $\overline{\mu}_0\in\Pi_{j=1}^N\Delta_w(S^j)$ the multi-population discrete-time mean-field game with total payoff defined with $r^i$, $Q^i$, $S^i$ and $A^i$, $i=1,\ldots,N$, has a Markov mean-field equilibrium.
\end{theorem}

{\em Proof}: Let us fix $\overline{\mu}_0$ and $M$ satisfying (A1''). 
Recall the notation used in the proof of Theorem \ref{thm:strong_discounted_mmfe_thm}
$$\Xi^i=\Pi_{t=0}^\infty \Delta_w^{(t)}(D^i),\quad
\Xi=\Pi_{i=1}^N\Xi^i.$$

Next, for any flow of measure-vectors $(\overline{\tau}):=(\overline{\tau}_0,\overline{\tau}_1,\ldots)\in\Xi$ and any $i\in\{ 1,\ldots,N\}$ let us define
$$V^{i,t}_{*(\overline{\tau})}(s):=\max_{\pi^i\in\mathcal{M}^i}\mathbb{E}^{\delta_s,\overline{Q_*},\pi^i}\sum_{k=t}^\infty r^i(s_k^i,a_k^i,\overline{\tau}_k),$$
that is, the optimal value at time $t\geq 0$ for the total-reward Markov decision process of player from population $i$ when the behaviour of all the other players at each stage is described by the flow $(\overline{\tau})$. Using the standard method of transforming a nonhomogeneous Markov decision process into a homogeneous one and Theorem 12 in \cite{FeinbergHuang}\footnote{Here we apply this theorem for weight function $V(s,t):=w(s)\alpha^{-t}$, $K:=\widetilde{L}_i$ and $\overline{c}=R$.} we may show that $V^i_{*(\overline{\tau})}(s)$ can be obtained from the optimal reward in the discounted Markov decision process with state space $S^i\times\mathbb{N}$ and:
\begin{enumerate}[(a)]
\item The function $\widetilde{\zeta}^i:S^i\times\mathbb{N}\times\Xi\rightarrow[0,\infty)$
$$\widetilde{\zeta}^i(s,t,(\overline{\tau})):=\sup_{\pi^i\in\mathcal{M}^i}\left[\sum_{k=0}^\infty\int_{S^i\setminus\{ s^*\}}w(s')\alpha^{-k}\left( Q^i_*\right)^k(ds'\mid s,\pi^i,(\overline{\tau}))\right]
$$
satisfying for $s\in S^i$, $t\geq 0$ and $(\overline{\tau})\in\Xi$, $\widetilde{\zeta}^i(s,t,(\overline{\tau}))\in [\alpha^{-t}w(s),\widetilde{L}_i\alpha^{-t}w(s)]$ with\footnote{Such a value will exist by (A2'') (b) and (A4'').} 
$$\widetilde{L}_i:=\sup_{(\overline{\tau})\in\Xi}\left\|\sum_{t=0}^\infty\int_{S^i\setminus\{ s^*\}}w(s')\alpha^{-t}\left( Q^i_*\right)^t(ds'\mid s,\pi^i,(\overline{\tau}))\right\|_w$$
\item The discount factor $\widetilde{\beta}:=\max_{j=1,\ldots,N}\frac{\widetilde{L}_j-1}{\widetilde{L}_j}$,
\item Modified one-stage rewards $\widetilde{r}^i_*$:
$$\widetilde{r}^i_*(s,t,a,(\overline{\tau})):=\frac{r^i(s,a,\overline{\tau}_t)}{\widetilde{\zeta}^i(s,t,(\overline{\tau}))},$$
\item Modified transition probabilities $\widetilde{Q}^i_{**}$ given by
\begin{align*}
&\widetilde{Q}^i_{**}(B\times\{ t'\}\mid s,t,a,(\overline{\tau})):=\\
&\left\{ \begin{array}{l}
\frac{1}{\widetilde{\beta}\widetilde{\zeta}^i(s,t,(\overline{\tau}))}\int_{B}\widetilde{\zeta}^i(s',t,(\overline{\tau}))Q^i_*(ds'\mid s,a,\overline{\tau}_t),\\
\quad\mbox{ if }B\in\mathcal{B}(S^i\setminus\{ s^*\}),(s,a)\in D^i\setminus\{ (s^*,a^*)\},t'=t+1,\\
1-\frac{1}{\widetilde{\beta}\widetilde{\zeta}^i(s,t,(\overline{\tau}))}\int_{B}\widetilde{\zeta}^i(s',t,(\overline{\tau}))Q^i_*(ds'\mid s,a,\overline{\tau}_t),\\
\quad\mbox{ if }B=\{ s^*\}),(s,a)\in D^i\setminus\{ (s^*,a^*)\},t'=t+1,\\
1,\mbox{ if }B=\{ s^*\}, (s,a)=(s^*,a^*),t'=t+1,\\
0,\mbox{ if }t'\neq t+1
\end{array}\right.
\end{align*}
\item Sets of feasible actions given by $\widetilde{A}^i(\cdot,t):=A^i(\cdot)$.
\end{enumerate}
In fact, if $V^i_{\widetilde{\beta}(\overline{\tau})}(s,t)$ denotes the optimal value in the modified (discounted) model, $$V^{i,t}_{*(\overline{\tau})}(s)=V^i_{\widetilde{\beta}(\overline{\tau})}(s,t)\widetilde{\zeta}^i(s,t,(\overline{\tau}))$$ 
for any $s\in S^i$ and $t\geq 0$. Moreover, optimal stationary strategies in the new model (which exist by Theorem 12 in \cite{FeinbergHuang}) correspond to optimal Markov strategies in the original one. Finally, repeating the arguments used in the proof of Lemma \ref{lem:total_V_cont} (the assumptions (A1--A4) used there are satisfied with $w(s)$ replaced by $\widetilde{w}(s,t):=w(s)\alpha^{-t}$, which clearly is a moment function on $S\times\mathbb{N}$) we can show that $V^i_{\widetilde{\beta}\cdot}(\cdot,t)$ and $\widetilde{\zeta}^i(\cdot,t,\cdot)$ are continuous and their $\widetilde{w}$-norms are bounded by $\widetilde{L}:=\frac{R\max_{j=1,\ldots,N}\widetilde{L}_j}{1-\widetilde{\beta}}$, which implies that for any $s\in S^i$ and $t\geq 0$, $V^{i,t}_{*(\overline{\tau})}(s)\leq \widetilde{L}\lambda^tw(s)$.

Next we define the correspondence $\widetilde{\Psi}^i_*$ from $\Xi$ into $\Xi^i$, ($i=1,\ldots,N$). In order to do it, we recall the definition of correspondence $\widetilde{\Theta}^i$ introduced in the proof of Theorem \ref{thm:strong_discounted_mmfe_thm}
$$
\widetilde{\Theta}^i((\overline{\tau}))=\left\{ (\eta^i)\in\Xi^i\! : \left(\eta^i_0\right)_{S^i}\!=\!\mu_0^i\mbox{ and }\!\left(\eta^i_t\right)_{S^i}(\cdot)\!=\!\int_{D^i}\! Q^i(\cdot\mid s,a,\overline{\tau}_{t-1})\tau^i_{t-1}(ds\times da)
\right\}\!.
$$
Now let
\begin{align*}
&\widetilde{\Psi}^i_*((\overline{\tau})):=
\left\{ (\eta^i)\in\widetilde{\Theta}^i((\overline{\tau})): \int_{S^i}V^{i,t-1}_{*(\overline{\tau})}(s)\left(\eta^i_t\right)_{S^i}(ds)\right.\\
&=\left.\int_{D^i}\left[ r^i(s,a,\overline{\tau}_{t-1})+\int_{S^i}V^{i,t}_{*(\overline{\tau})}(s')Q^i_*(ds'\mid s,a,\overline{\tau}_{t-1})\right]\eta^i_{t-1}(ds\times da)\mbox{ for }t\geq 1
\right\}
\end{align*}
As we have shown in the proof of Theorem \ref{thm:strong_discounted_mmfe_thm}\footnote{The assumptions (A1') and (A2') used there are in fact slightly stronger than (A1'') and (A2'') assumed here, so the proof remains valid in the present case.}, for each $i$, the correspondence $\widetilde{\Theta}^i$ has non-empty values and closed graph.
We next prove that for any fixed $i\in\{ 1,\ldots,N\}$, $\widetilde{\Psi}^i_*$ has similar properties. 

We start by noting, that if $\pi^i_*$ is an optimal deterministic Markov policy in the optimization problem of a player from population $i$ maximizing his total reward when the behaviour of all the other players at each stage is described by the flow $(\overline{\tau})$, then for each $t\geq 1$ and $s\in S^i$ it satisfies
$$V^{i,t-1}_{*(\overline{\tau})}(s)=\left[ r^i(s,\pi^i_{*t-1}(s),\overline{\tau}_{t-1})+\int_{S^i}V^{i,t}_{*(\overline{\tau})}(s')Q^i_*(ds'\mid s,\pi^i_{*t-1}(s),\overline{\tau}_{t-1})\right],$$
which implies that for any $t\geq 1$,
\begin{eqnarray*}
&&\int_{S^i}V^{i,t-1}_{*(\overline{\tau})}(s)\left( \Pi^i_{t-1}(\pi^i_*,\Phi^i(\cdot\mid \overline{\tau}_{t-1}))\right)_{S^i}(ds)=\int_{S^i}V^{i,t-1}_{*(\overline{\tau})}(s)\Phi^i(\cdot\mid \overline{\tau}_{t-1})(ds)\\
&=&\int_{S^i}\left[ r^i(s,\pi^i_{*t-1}(s),\overline{\tau}_{t-1})+\int_{S^i}V^{i,t}_{*(\overline{\tau})}(s')Q^i_*(ds'\mid s,\pi^i_{*t-1}(s),\overline{\tau}_{t-1})\right] \Phi^i(\cdot\mid \overline{\tau}_{t-1})(ds)\\
&=&\int_{D^i}\left[ r^i(s,a,\overline{\tau}_{t-1})+\int_{S^i}V^{i,t}_{*(\overline{\tau})}(s')Q^i_*(ds'\mid s,a,\overline{\tau}_{t-1})\right]\Pi^i_{t-1}(\pi^i_*,\Phi^i(\cdot\mid \overline{\tau}_{t-1}))(ds\times da),
\end{eqnarray*}
Hence, measure flow $(\eta^i)$ defined by (\ref{eq:rec}) with $\pi^i:=\pi^i_*$ is an element of $\widetilde{\Psi}^i_*((\overline{\tau}))$.

In the penultimate part of the proof we show that the graph of $\widetilde{\Psi}^i_*$ is closed.
Let $\left\{(\overline{\tau}^{(n)})\right\}_{n\geq 0}$ in $\Xi$ and $\left\{(\eta^{i,(n)})\right\}_{n\geq 0}$ in $\Xi^i$ be convergent sequences such  that $(\eta^{i,(n)})\in\widetilde{\Psi}^i_*((\overline{\tau}^{(n)}))$ for each $n$. Moreover, let $(\eta^{i,(n)})\Rightarrow(\eta^i)$ and $(\overline{\tau}^{(n)})\Rightarrow(\overline{\tau})$ as $n\rightarrow\infty$ for some $(\eta^i)\in\Xi^i$ and $(\overline{\tau})\in\Xi$. Since the graph of $\widetilde{\Theta}^i$ is closed, showing that  $\widetilde{\Psi}^i_*$ has the same property only requires proving that the equalities defining $\widetilde{\Psi}^i_*$ hold for $(\eta^i)$ and $(\overline{\tau})$.
Let us fix $t\geq 1$.
From the definition of $\widetilde{\Psi}^i_*$ we know that for each $n$
\begin{eqnarray}
\label{eq:TBellman_n}
&&\int_{S^i}V^{i,t-1}_{*(\overline{\tau}^{(n)})}(s)\left(\eta^{i,(n)}_t\right)_{S^i}(ds)\\
&=&\int_{D^i}\left[ r^i(s,a,\overline{\tau}^{(n)}_{t-1})+\int_{S^i}V^{i,t}_{*(\overline{\tau}^{(n)})}(s')Q^i_*(ds'\mid s,a,\overline{\tau}^{(n)}_{t-1})\right]\eta^{i,(n)}_{t-1}(ds\times da)\nonumber
\end{eqnarray}
Then by the continuity of $Q^i_*$ and $V^{i,t}_{*(\overline{\tau}^{(n)})}$ and Theorem 3.3 in \cite{Serfozo} (see also Remark 3.4 (ii) there -- by (A2'') (b) the assumption presented there is true for $g=\widetilde{L}\alpha^tw$), $\int_{S^i}V^{i,t}_{*(\overline{\tau}^{(n)})}(s')Q^i_*(ds'\mid \cdot,\cdot,\overline{\tau}^{(n)}_{t-1})$ converges continuously to $\int_{S^i}V^{i,t}_{*(\overline{\tau}}(s')Q^i_*(ds'\mid \cdot,\cdot,\overline{\tau}_{t-1})$. As also $r^i(\cdot,\cdot,\overline{\tau}^{(n)}_{t-1})$ converges continuously to $r^i(\cdot,\cdot,\overline{\tau}_{t-1})$ by (A1''), using Theorem 3.3 in \cite{Serfozo} again (again with $g=\widetilde{L}\alpha^tw$, which satisfies the assumption given in Remark 3.4 (ii) by (A1'') and (A2'') (b)), we can pass to the limit in (\ref{eq:TBellman_n}), obtaining
\begin{eqnarray*}
&&\int_{S^i}V^{i,t-1}_{*(\overline{\tau})}(s)\left(\eta^i_t\right)_{S^i}(ds)\\
&=&\int_{D^i}\left[ r^i(s,a,\overline{\tau}_{t-1})+\int_{S^i}V^{i,t}_{*(\overline{\tau})}(s')Q^i_*(ds'\mid s,a,\overline{\tau}_{t-1})\right]\eta^i_{t-1}(ds\times da)
\end{eqnarray*}
As $t$ was arbitrary, this ends the proof that the graph of $\widetilde{\Psi}^i_*$ is closed.

The remainder of the proof is identical to the argument presented in the proof of Theorem \ref{thm:strong_discounted_mmfe_thm}:
%
We define the correspondence from $\Xi$ into itself:
$$\widetilde{\Psi}_*((\overline{\tau})):=\widetilde{\Psi}^1_*((\overline{\tau}))\times\ldots\times \widetilde{\Psi}^N_*((\overline{\tau})).$$
By what we have shown, $\widetilde{\Psi}_*$ has nonempty values and its graph is closed. Convexity of values of $\widetilde{\Psi}_*$ is obvious. As we know, 
$\Xi$ is compact in product topology. Therefore,
Glickberg's fixed point theorem \cite{Glicksberg52} implies that $\widetilde{\Psi}_*$ has a fixed point. Suppose $(\overline{\tau^*})$ is this fixed point. Disintegrating $\tau^{*i}_t$ gives for $i=1,\ldots,N$ and $t=0,1,\ldots$ stochastic kernels $\pi^{*i}_t$ and measures $\mu^{*i}_t$ which (after similar modifications as in the proof of Theorem \ref{thm:strong_discounted_smfe_thm}) correspond to the Markov strategies and global state flows in mean-field Markov equilibrium in the total-reward game.
$\Box$


\begin{remark}
It should be noted here that Theorems \ref{thm:strong_total_smfe_thm} and \ref{thm:strong_total_mmfe_thm} applied to a game with a single population extend the existing results about the existence of equilibria in single-population total-reward mean-field games. The only results of such type in the literature appear in \cite{PiWEA} and concern model with finite state and action spaces. In \cite{PiWEA} it was assumed that the probability of reaching $s^*$ within some fixed number of stages is for any strategies used by the players bigger than some $p_0>0$. Assumptions (A4) and (A4'') used here can be seen as counterparts of this assumption for our model with weight function $w$ applied to the states. It is easy to see that in the finite state and action case (A4) reduces to the simpler assumption described above.
\end{remark}

\begin{remark}
It is also worth noting that some total-reward mean-field game models that are not directly considered in the article can be treated as specific cases of our framework. Firstly, a total reward game without a replacement of dead players by new-born ones is a specific case with $s^*$ being an absorbing state. For such a case the existence of Markov mean-field equilibrium is guaranteed by Theorem \ref{thm:strong_total_mmfe_thm} without any modification of assumptions.

Another case that can be treated as an instance of our model is the total reward game with a finite horizon. In this case the application of our framework requires replacing individual state spaces of players of each population $S^i$ with $\left( (S^i\setminus\{ s^*\})\times \{ 0,\ldots,T\}\right)\cup\{ s^*\}$ (with $T$ denoting the time horizon). If the stochastic kernel $Q^i$ denotes the transition probability for population $i$ in the original model, the transitions in the modified one $\widehat{Q}^i$ are defined completely with the following formulas:
$$\widehat{Q}^i(B\times \{ t+1\}\mid (s,t),a,\overline{\tau})=Q^i(B\mid s,a,\overline{\tau}_S)$$
for $s\in S^i$, $a\in A^i(s)$, $t<T$ and $B\in\mathcal{B}(S^i)$, with $\overline{\tau}_S=(\tau^1_{S^1},\ldots,\tau^N_{S^N})$ being a vector of marginals of measures $\tau^i$ on original individual state spaces $S^i$, and
$$\widehat{Q}^i(\{ s^*\}\mid (\cdot,T),\cdot,\cdot)\equiv 1\equiv \widehat{Q}^i(\{ s^*\}\mid s^*,a^*,\cdot).$$
Then the single-stage rewards are defined independently of the time components of both the individual and the global state. The existence of Markov mean-field equilibria is then assured by Theorem \ref{thm:strong_total_mmfe_thm} under assumptions (A1)--(A3) on original primitives of the model (which correspond to (A1''), (A2'') and (A3) for the modified one). (A4'') is satisfied automatically.
A similar transformation allows for considering the non-stationary model with finite horizon as well as the case when the time horizon of each individual is a random variable with a finite expected value independent from the Markov chain of his individual states.
\end{remark}

\section{Concluding remarks.}

In the paper we have presented a model of discrete-time mean-field game with several populations of players. 
Games of this type have been studied in the literature in the discrete-time setting. The main results presented in this article are stationary and Markov mean-field equilibrium existence theorems for two payoff criteria: $\beta$-discounted payoff and total payoff proved under some rather general assumptions on one-step reward functions and individual transition kernels of the players. 
It is also worth noting that the games with total payoff have only been studied in finite state space case, hence, the results presented here also extend those for total-payoff mean-field games with a single population. 
The article is the first of two papers on multiple-population discrete-time mean-field games with discounted or total payoff. In the second one we provide theorems showing that under some additional assumptions equilibria obtained in the mean-field models are approximate equilibria in their $n$-person counterparts when $n$ is large enough. We also plan further research on the topic of discrete-time mean-field games with multiple populations of players which will concentrate on games with long-run average reward.


\section*{Declarations}
\begin{itemize}
\item {\bf Funding}
This work was supported by the NCN Grant no 2016/23/B/ST1/00425.
\item {\bf Conflict of interest}
The author declares that he has no conflict of interest.
\end{itemize}

\end{document}